\documentclass[10pt]{amsart}

\usepackage{amsmath, amssymb, latexsym, enumerate}
\input{diagrams}

\newtheorem{thm}[equation]{Theorem}
\newtheorem{pro}[equation]{Proposition}
\newtheorem{cor}[equation]{Corollary}
\newtheorem{lem}[equation]{Lemma}
\newtheorem{dfn}[equation]{Definition}

\makeatletter\@addtoreset{equation}{section}\makeatother

\theoremstyle{definition}
\newtheorem{exa}[equation]{Example}

\newcommand{\sands}{\quad\text{and}\quad}
\newcommand{\spandsp}{\qquad\text{and}\qquad}
\newcommand{\co}{\colon}
\newcommand{\coo}{\,\colon\,}
\newcommand{\rta}{\rightarrow}
\renewcommand{\o}{\otimes}
\renewcommand{\=}{\,=\,}

\newcommand{\isom}{\cong}
\newcommand{\pr}[2]{\ensuremath{\langle {#1,#2}\rangle}}
\newcommand{\sub}[1]{_{\text{${\scriptscriptstyle #1}$}}}

\newcommand{\K}{\ensuremath{K}}   

\newcommand{\as}{\ensuremath{\text{$A\sub S$}}}
\newcommand{\at}{\ensuremath{\text{$A\sub T$}}}
\newcommand{\au}{\ensuremath{\text{$A\sub U$}}}

\newcommand{\bs}{\ensuremath{\text{$B\sub S$}}}
\newcommand{\bt}{\ensuremath{\text{$B\sub T$}}}
\newcommand{\cs}{\ensuremath{\text{$C\sub S$}}}
\newcommand{\ct}{\ensuremath{\text{$C\sub T$}}}

\newcommand{\us}{\ensuremath{\text{$U\sub S$}}}
\newcommand{\ut}{\ensuremath{\text{$U\sub T$}}}
\newcommand{\vs}{\ensuremath{\text{$V\sub S$}}}
\newcommand{\vt}{\ensuremath{\text{$V\sub T$}}}

\newcommand{\rkm}[1]{\ensuremath{\rho_{\scriptscriptstyle M}(#1)}}
\newcommand{\rkmd}[1]{\ensuremath{\rho_{\scriptscriptstyle M^*}(#1)}}
\newcommand{\rkn}[1]{\ensuremath{\rho_{\scriptscriptstyle N}(#1)}}
\newcommand{\rkl}[1]{\ensuremath{\rho_{\scriptscriptstyle L}(#1)}}
\newcommand{\nlm}[1]{\ensuremath{\nu_{\scriptscriptstyle M}(#1)}}
\newcommand{\nlp}[1]{\ensuremath{\nu_{\scriptscriptstyle P}(#1)}}
\newcommand{\nlmd}[1]{\ensuremath{\nu_{\scriptscriptstyle M^*}(#1)}}
\newcommand{\nln}[1]{\ensuremath{\nu_{\scriptscriptstyle N}(#1)}}
\newcommand{\rlm}[1]{\ensuremath{\lambda_{\scriptscriptstyle M}(#1)}}

\newcommand{\rln}[1]{\ensuremath{\lambda_{\scriptscriptstyle N}(#1)}}
\newcommand{\nll}[1]{\ensuremath{\nu_{\scriptscriptstyle L}(#1)}}
\newcommand{\clm}[1]{\ensuremath{c\ell_{\scriptscriptstyle M}(#1)}}
\newcommand{\clq}[1]{\ensuremath{c\ell_{\scriptscriptstyle Q}(#1)}}
\newcommand{\cln}[1]{\ensuremath{c\ell_{\scriptscriptstyle N}(#1)}}
\newcommand{\cll}[1]{\ensuremath{c\ell_{\scriptscriptstyle L}(#1)}}

\newcommand{\mcat}{\ensuremath{\mathbf M}}
\newcommand{\idt}[1]{\ensuremath{id\sub{#1}}}

\newcommand{\fs}{\ensuremath{\mathcal F}}
\newcommand{\dlat}{\ensuremath{\mathcal D}}

\newcommand{\frp}{\mathbin{\Box}}

\newcommand{\ca}{\ensuremath{\mathcal A}}
\newcommand{\cb}{\ensuremath{\mathcal B}}

\newcommand{\cf}{\ensuremath{\mathcal F}}

\newcommand{\cm}{\ensuremath{\mathcal M}}
\newcommand{\cn}{\ensuremath{\mathcal N}}
\newcommand{\cu}{\ensuremath{\mathcal U}}
\newcommand{\calt}{\ensuremath{\mathcal T}}

\newcommand{\frmod}[1]{\ensuremath{K\{{#1}\}}}     
\newcommand{\sect}[3]{\ensuremath{\binom{#1}{#2,#3}}}  
\newcommand{\msect}[3]{\ensuremath{\binom{#1}{#2,\dots,#3}}} 
\newcommand{\dual}{\ensuremath{^*}}
\newcommand{\pt}{\ensuremath{I}}            
\newcommand{\lp}{\ensuremath{Z}}

\newcommand{\un}[2]{\ensuremath{U_{#1,#2}}} 
\newcommand{\ist}[1]{\ensuremath{\text{\rm Isth}(#1)}}
\newcommand{\loo}[1]{\ensuremath{\text{\rm Loop}(#1)}}
\newcommand{\lift}[1]{\ensuremath{\text{\rm L$^{#1}$}}}
\newcommand{\li}{\ensuremath{\text{\rm L}}}
\newcommand{\trun}[1]{\ensuremath{\text{\rm T$^{#1}$}}}
\newcommand{\tr}{\ensuremath{\text{\rm T}}}

\begin{document}

\title[A unique factorization theorem for matroids]
{A unique factorization theorem for matroids} 
\author{Henry Crapo}
\author{William Schmitt}
\thanks{Schmitt partially supported by NSA grant 02G-134}
\email{crapo@ehess.fr and wschmitt@gwu.edu}
\begin{abstract}
   We study the combinatorial, algebraic and geometric properties
of the free product operation on matroids.  After giving cryptomorphic 
definitions of free product in terms of independent sets, bases,
circuits, closure, flats and rank function, we show that free
product, which is a noncommutative operation, is associative 
and respects matroid duality.  The free product of matroids $M$ and $N$
is maximal with respect to the weak order among matroids having $M$
as a submatroid, with complementary contraction equal to $N$.  
Any minor of the free product of $M$ and $N$ is a free product
of a repeated truncation of the corresponding minor of $M$ with
a repeated Higgs lift of the corresponding minor of $N$.  
We characterize, in terms of their cyclic flats, matroids that are 
irreducible with respect to free product, and prove that the factorization
of a matroid into a free product of irreducibles is unique up to 
isomorphism.   We use these results to determine, for \K\ a field
of characteristic zero, the structure of
the minor coalgebra $\frmod\cm$ of a family of matroids \cm\ that
is closed under formation of minors and free products: namely,
\frmod\cm\ is cofree, cogenerated by the set of irreducible matroids
belonging to \cm. 
\end{abstract}

\maketitle
\begin{center}
  {\it For Denis Higgs, who gave us the lift}.
\end{center}
\vspace{.25in}
\section{Introduction}

We introduced the free product of matroids in a short article
(\cite{crsc:fpm}), in which we used it to settle the conjecture by
Welsh (\cite{we:bnm}) that $f_{n+m}\geq f_n\cdot f_m$, where $f_n$ is
the number of distinct isomorphism classes of matroids on an
$n$-element set.  Free product is, in a categorical sense, dual to the
direct sum operation, and has properties that are in striking contrast
to those of other, better known, binary operations on matroids; most
significantly, it is noncommutative. In the present article we
initiate a systematic study of the combinatorial, algebraic and
geometric properties of this new operation.  Our main results include
a characterization, in terms of cyclic flats, of matroids that are
irreducible with respect to free product, and a unique factorization
theorem: every matroid factors uniquely, up to isomorphism, as a free
product of irreducible matroids.  Hence the set
of all isomorphism classes of matroids, equipped with the binary
operation induced by free product, is a free monoid, generated by the
isomorphism classes of irreducible matroids.

Although we first defined the free product as such in \cite{crsc:fpm},
we first became aware of it earlier, while investigating, in
\cite{crsc:fsa}, the {\it minor coalgebra} of a minor-closed family of
matroids.  This coalgebra has as basis the set of all isomorphism
classes of matroids in the given family, with coproduct of a matroid
$M=M(S)$ given by $\sum\sub{A\subseteq S}M|A\o M/A$, where $M|A$ is
the submatroid obtained by restriction to $A$ and $M/A$ is the
complementary contraction.  If the family is also closed under
formation of direct sums then its minor coalgebra is a Hopf algebra,
with product determined on the basis of matroids by direct sum.  These
Hopf algebras, and analogous Hopf algebras based on families of graphs,
were introduced in \cite{sc:iha}, as examples of the
more general construction of incidence Hopf algebra.  In the dual of
the minor coalgebra, the {\it minor algebra},
the product of matroids $M$ and $N$ (dual basis
elements) is a linear combination of those matroids having some
restriction isomorphic to $M$, with complementary contraction
isomorphic to $N$; the coefficient of $L=L(U)$ being the number of
subsets $A\subseteq U$ such that $L|A\isom M$ and $L/A\isom N$.  In
the weak map order, the set of matroids appearing with nonzero
coefficient in this product has a minimum element, given by the direct
sum $M\oplus N$, and also has a maximum element, which we denote by
$M\frp N$; this is the free product of
$M$ and $N$.

After discussing a few preliminaries in the following short section, we
begin Section \ref{sec:crypt} by recalling from \cite{crsc:fpm} 
the definition, in terms of independent sets, of the free product.
As a next step, dictated by the culture of matroid theory, we give
cryptomorphic definitions of the free product in terms of bases, circuits,
closure, flats and rank function.  These various characterizations
allow us to demonstrate, in Sections \ref{sec:basic} and \ref{sec:minors}, 
a number of fundamental properties of free product. In particular:
free product satisfies the extremal property mentioned above, that is,
$M\frp N$ is maximal in the weak order among
matroids having a submatroid equal to $M$, with complementary
contraction equal to $N$;  free product
is associative, and commutes with matroid duality; and 
any minor of a free product $M\frp N$ is itself a
free product, namely, the free product of a repeated truncation of a
minor of $M$ with a repeated Higgs lift of a minor of $N$.

We begin Section \ref{sec:uft} by giving a characterization of the
cyclic flats of a free product, and making the key definition of {\it
  free separator} of a matroid $M(S)$ as a subset of $S$ that is
comparable by inclusion to all cyclic flats of $M$.  We then prove the
theorem that $M$ factors as a free product $P(U)\frp Q(V)$ if and only
if the set $U$ is a free separator of $M$.  As a consequence, we find
that a nonuniform matroid $M(S)$ is irreducible if and only if the
complete sublattice $\dlat (M)$ of the Boolean algebra $2^S$ generated
by the cyclic flats of $M$ has no {\it pinchpoint}, that is,
single-element crosscut, other than $\emptyset$ and $S$. (Uniform matroids
factor completely, into single-element matroids.) In order to
examine free product factorization of matroids in detail, we turn our
attention to the set $\fs (M)$ of all free separators of a matroid
$M(S)$, which, partially ordered by inclusion, is also a sublattice of
$2^S$.  By the theorem mentioned above, there is a one-to-one
correspondence between chains from $\emptyset$ to $S$ in $\fs (M)$ and
factorizations $M=M_1\frp\cdots\frp M_k$, according to which $M_i$ is
the minor of $M$ determined by the $i${\raisebox{.3ex}
  {$\scriptstyle\rm th$}} interval in the corresponding chain.
Factorizations of $M$ into irreducibles thus correspond to maximal
chains in $\fs (M)$.

We define the {\it primary flag} $\calt\sub M$ of a matroid $M$ as the
chain $T_0\subset\cdots\subset T_k$ of pinchpoints in the
lattice $\dlat (M)$.  We show that $\calt\sub M$ is also the chain of
pinchpoints in $\fs (M)$ and, furthermore, that the
intersection of the lattices $\fs (M)$ and $\dlat (M)$ is precisely
$\calt\sub M$.  These results, together with a proposition
characterizing the intervals $[T_{i-1}, T_i]$ in $\fs (M)$, allow us
to prove that the free product factorization $M=M_1\frp\cdots\frp M_k$
corresponding to the chain $\calt\sub M$ is the unique
factorization of $M$ having the property that each $M_i$ is either
irreducible, or maximally uniform (in the sense that no free product
of consecutive $M_i$'s is uniform).  From this fundamental result, our
main theorem quickly follows: every matroid factors uniquely up to
isomorphism as a free product of irreducible matroids.

In Section \ref{sec:algebra}, we use the
unique factorization theorem, together with the extremal property of
free product with respect to the weak order, to show that
for any class \cm\ of matroids closed under the formation of minors
and free products, the minor coalgebra of \cm\ is cofree, cogenerated by
the isomorphism classes of irreducible matroids in \cm.
Any minor-closed class of matroids defined by the exclusion of a set
of irreducible minors will therefore 
generate a minor coalgebra that is
cofree.  This is not the case for certain well-studied classes such as
binary or unimodular matroids, because the four point line factors (as
the free product of four one-element matroids).  But for an infinite
field $F$ the class of $F$-representable matroids is closed under free
product and hence its minor coalgebra is cofree.

In conclusion, we sketch in Section \ref{sec:twist} a development
whereby the minor coalgebra of a free product and minor-closed family
of matroids forms a (self-dual) Hopf algebra in an appropriate braided
monoidal category.
\section{preliminaries}
We denote the disjoint union of sets $S$ and $T$ by $S+T$, the set
difference by $S\backslash T$, and the intersection $S\cap T$ by either
$S\sub T$ or $T\sub S$.  If $T$ is a singleton set $\{a\}$, we write
$S+a$ and $S\backslash a$, respectively for $S+T$ and $S\backslash T$.
We write $M=M(S)$ to indicate that $M$ is a matroid with ground set $S$;
in the case that $S=\{a\}$ is a singleton set
we write $M(a)$ instead of $M(S)$.
We denote the rank and nullity functions of $M$ by $\rho\sub M$ and
$\nu\sub M$, respectively, and denote 
by $\lambda\sub M$ the {\it rank-lack} function on
$M$, given by $\rlm A= \rho (M) - \rkm A$, for all $A\subseteq S$,
where $\rho (M)=\rkm S$ is the rank of $M$.

Given a matroid $M(S)$ and $A\subseteq S$, we write $M|A$ for the
restriction of $M$ to $A$, that is, the matroid on $A$ obtained by
deleting $S\backslash A$ from $M$, and we write $M/A$ for the matroid
on $S\backslash A$ obtained by contracting $A$ from $M$.  For all
$A\subseteq B\subseteq S$, we denote the minor
$(M|B)/A=(M/A)|(B\backslash A)$ by $M(A,B)$.

 For any set $S$, the {\it free matroid} $I(S)$ and the
{\it zero matroid} $Z(S)$ are, respectively, the unique matroids on
$S$ having nullity zero and rank zero.  In other words, if $|S|=n$,
then $I(S)$ is the uniform matroid $U_{n,n}(S)$ and $Z(S)$ is
the uniform matroid $U_{0,n}(S)$.  
We refer the reader to Oxley \cite{ox:mt} and Welsh \cite{we:mt}
for any background on matroid theory that might be needed.

\section{The free product: cryptomorphic definitions}\label{sec:crypt}

\begin{dfn}[\cite{crsc:fpm}]
{\rm
The {\it free product} of matroids $M(S)$ and $N(T)$ is
the matroid $M\frp N$ defined on the set $S+T$ whose collection of
independent sets is given by}
$$
\{ A\subseteq S+T\coo\text{\rm $A\sub S$ is independent in $M$
and $\rlm\as \geq\nln\at$}\}.
$$
\end{dfn}
The first two propositions of \cite{crsc:fpm} show that
$M\frp N$ is indeed a matroid, which contains $M$ and $N$
as complementary minors; specifically, if the ground set of $M$
is $S$, then
\begin{equation}
  \label{eq:univ1}
  (M\frp N)|S = M \spandsp\ (M\frp N)/S = N.
\end{equation}

\begin{pro}\label{pro:bases}
The collection of bases of $M(S)\frp N(T)$ is given by
$$
\{ A\subseteq S+T\coo\text{\rm \as\ is independent in $M$, \at\ spans
  $N$, and $\rlm\as=\nln\at$}\}.
$$
\end{pro}
\begin{proof}
  The result follows directly from the definition of
the free product.
\end{proof}
Note that it follows immediately from the characterization of the
bases of $M\frp N$ that $\rho (M\frp N)=\rho (M)+\rho (N)$,
for all $M$ and $N$.
\begin{exa}\label{exa:frps}
  Let $S=\{e,f,g\}$ and $T=\{a,b,c,d\}$, and suppose that $M(S)$ is
a three-point line, and $N(T)$ consists of two double points
$ab$ and $cd$.  The free products $I(e)\frp N(T)$ and $M(S)\frp N(T)$
are shown below:

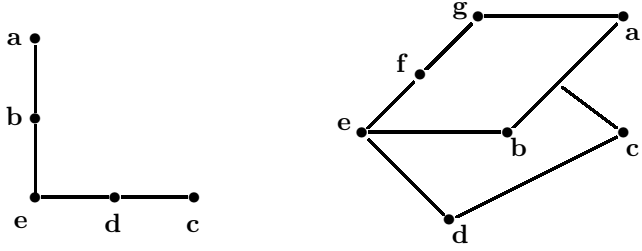
\newcommand{\hllab}[1]
{\hspace{3.3ex}\text{\raisebox{0ex}{$\bf #1$}}}
\newcommand{\eelab}[1]
{\hspace{2.ex}\text{\raisebox{.8ex}{$\bf #1$}}}
\newcommand{\elab}[1]
{\hspace{-3.1ex}\text{\raisebox{-1.5ex}{$\bf #1$}}}
\vspace{1.5ex}
\begin{center}
\mbox{
\begin{diagram}[PostScript=Rokicki,abut,height=1.5em,width=1.5em,tight,thick]
&\hllab{a}&\bullet\\
&&\dLine\\
&\hllab{b}&\bullet\\
&&\dLine\\
&{\hspace{4.45ex}\text{\raisebox{-2.1ex}{$\bf e$}}}
&\bullet&\rLine&\bullet&\rLine
&\bullet&
{\hspace{-14.3ex}\text{\raisebox{-2.5ex}{$\bf d$\hspace{5.7ex}$\bf c$}}}
\end{diagram}}\qquad\qquad
\mbox{
\begin{diagram}[PostScript=Rokicki,abut,height=1.1em,width=1.1em,tight,thick,balance]
&&&&\eelab{g}&\bullet&&\rLine&&&\bullet&\elab{\!a}\\
&&&&\ruLine(2,2)&&&&&\ruLine(4,4)\\
&&\eelab{f}&\bullet&&&&
\rdLine[leftshortfall=2.6em](4,3)\\
&&\ruLine(2,2)\\
\eelab{e}&\bullet&\rLine&&&&\bullet&\elab{b}&&&
\bullet&\elab{\!c}\\
&&\rdLine(3,3)&&&&&&&\ldLine(6,3)\\\\
&&&&\bullet&\elab{d}
\end{diagram}}
\end{center}
\vspace{1.5ex}
According to Proposition \ref{pro:bases}, the matroid $I\frp N$ has as bases 
all three-element subsets of $\{a,b,c,d\}$, together with
all sets of the form $\{e,x,y\}$, where $x\in\{a,b\}$ and $y\in\{c,d\}$;
while the bases of $M\frp N$ are the sets of the form $A\cup B$, 
with $A\subseteq S$, $B\subseteq T$,
and either
\begin{enumerate}[(i)]
\item $A=\emptyset$\, and\, $B=T$,
\item $|A|=1$\, and\, $|B|=3$,\, or
\item $|A|=2$\, and\, $|B|=2$,\, with $B$ not equal to $\{a,b\}$ or 
$\{c,d\}$.
\end{enumerate}
\end{exa}

\begin{pro}\label{pro:rank}
The rank function of $L=M(S)\frp N(T)$ is given by
$$
\rkl A \= \rkm\as + \rkn\at +\min\{\rlm\as,\nln\at\}, 
$$
for all $A\subseteq S+T$. 
\end{pro}
\begin{proof}
Suppose that $A\subseteq S+T$ and that $\rlm\as\geq\nln\at$.
Then for any basis $B$ of $M|\as$, the set $B\cup\at$ is a basis
for $L|A$, and thus $\rkl A=|B\cup\at|=|B|+|\at|=\rkm\as +\rkn\at +\nln\at$.

If $\rlm\as\leq\nln\at$, choose $C\subseteq\at$ such that 
$\rkn{C}=\rkn\at$ and $\nln{C}=\rlm\as$ and note
that we then have $|C|=\rkn{C}+\nln{C}=
\rkn\at+\rlm\as$. If $B$ is a basis for $M|\as$, then
$B\cup C$ is a basis for $L|A$, and thus $\rkl A =
|B\cup C| = \rkm\as + \rkn\at+\rlm\as$.
\end{proof}

It follows immediately that  the 
nullity function of $L=M(S)\frp N(T)$ is given by 
\begin{equation}
  \label{eq:null}
\nll A \= \nlm\as + \nln\at -\min\{\rlm\as,\nln\at\}, 
\end{equation}
for all $A\subseteq S+T$, and similarly for the rank-lack
function.

\begin{pro}\label{pro:closure}
  The closure operator on $L=M(S)\frp N(T)$ is given by 
$$
\cll A \=
\begin{cases}
  \clm\as\cup\at , & \text{if $\rlm\as>\nln\at$},\\
  S\cup\cln\at , &\text{if $\rlm\as\leq\nln\at$},
\end{cases}
$$
for all $A\subseteq S+T$.
\end{pro}
\begin{proof}
  Suppose that $\rlm\as>\nln\at$. According to Proposition \ref{pro:rank},
the rank of $A$ in $L$ is given by $\rkl A = \rkm\as +|\at|$, and
if $B=A\cup x$, for any $x\in S+T$, then $\rlm\bs\geq\nln\bt$,
and we have $\rkl B = \rkm\bs +|\bt|$.  
Hence $x\in\cll A$ if and only if $\rkm\as +|\at| = \rkm\bs + |\bt|$,
that is, if and only if $x\in\clm\as\cup\at$.

Suppose that $\rlm\as\leq\nln\at$.  If $B=A\cup x$, for any
$x\in S+T$, then $\rlm\bs\leq\nln\bt$ and thus, by Proposition
\ref{pro:rank}, $\rkl A=\rho (M)+
\rkn\at$ and $\rkl B = \rho (M)+\rkn\bt$.  Hence $x\in\cll A$ if
and only if $\rkn\at =\rkn\bt$, that is, if and only if
$x\in S\cup\cln\at$.
\end{proof}
As a corollary, we obtain the following description of
the flats of a free product in terms of the flats of its factors.  
\begin{cor}\label{cor:flats}
 Suppose that $L=M(S)\frp N(T)$ and $A\subseteq S+T$.  If $\rlm\as>
\nln\at$, then $A$ is a flat of $L$ if and only if \as\ is a flat
of $M$;  if $\rlm\as\leq\nln\at$, then $A$ is a flat of $L$ if and only
if $\as=S$ and \at\ is a flat of $N$. 
\end{cor}

\begin{pro}\label{pro:circuit}
  A set $C\subseteq S+T$ is a circuit in $L=M(S)\frp N(T)$ if and only
if $C\subseteq S$ and
$C=\cs$ is a circuit in $M$, or \cs\ is independent in $M$,
the restriction $N|\ct$ is isthmusless, and $\rlm\cs +1=\nln\ct$.
\end{pro}
\begin{proof}
  By the definition of free product, a subset $C$ of $S+T$ is
dependent in $L$ if and only if \cs\ is dependent in $M$ or
$\rlm\cs<\nln\ct$.  A minimal set with this property is either a circuit
in $M$, or a minimal set with \cs\ independent in $M$ but with 
$\rlm\cs < \nln\ct$, that is, a set such that $\rlm\cs +1 =\nln\ct$.
If such a set $C$ were such that the restriction $N|\ct$ were to have an
isthmus $d$, then $C$ would not be minimal, since we would have
$\nln\ct = \nln{\ct\backslash d}$. 
\end{proof}

\section{Basic properties of the free product}\label{sec:basic}

We begin with a lemma showing that the asserted inequality between
$\rlm\as$ and $\nln\at$ in the definition of free product is in
fact a property of restrictions and complementary contractions in
arbitrary matroids.
\begin{lem}\label{lem:conineq}
   Given a matroid $L=L(S+T)$, let $M=L|S$ and $N=L/S$.  Then
$\rlm\as\geq\nln\at$, for all independent sets $A$ in $L$.
\end{lem}
\begin{proof}
  The rank function on the contraction $N=L/S$ is determined by
  $\rkn B = \rkl{B\cup S} -\rkl S= \rkl{B\cup S}- \rho (M)$,
for all $B\subseteq T$.  If $A\subseteq S+T$ is independent in $L$,
then $\rkl{\at\cup S}\geq |A|$, and so by the above formula,
$\rkn\at\geq |A|-\rho (M)$.  Thus we have
 $\nln\at = |\at|-\rkn\at\leq |\at| - (|A|-\rho (M))=\rlm\as$.
\end{proof}
By definition, the independent sets of the free product $M(S)\frp N(T)$ are
precisely those subsets of $S+T$ which, according to Lemma
\ref{lem:conineq}, are necessarily independent in any matroid
containing $M$ as a submatroid with complementary contraction $N$.
The following proposition expresses the consequent extremal, or
universal, property of the free product.
\begin{pro}\label{pro:univ}
  For any matroid $L=L(U)$, and $S\subseteq U$, the identity
map on $U$ is a rank-preserving weak map $L|S\frp L/S\rta L$.
\end{pro}
\begin{proof}
  Let $M=L|S$ and $N=N(T)=L/S$. If $A$ is independent in $L$, then
$\as$ is independent in $M$ and, by Lemma \ref{lem:conineq},
we have $\rlm\as\geq\nln\at$.  Hence $A$ is independent in $M\frp N$,
and so the identity map on $S+T$ is a weak map from $M\frp N$
to $L$, which is clearly rank-preserving.
\end{proof}

Roughly speaking, in a free product $L=M(S)\frp N(T)$, the submatroid
$L|T$ is the freest matroid, arranged in the most general
position possible relative to $M=L|S$ such that the contraction $L/S$
is equal to $N(T)$.  In the matroid $M(S)\frp N(T)$ of Example
\ref{exa:frps}, as long as $\{a,b\}$ and $\{c,d\}$ are each coplanar
with $S=\{e,f,g\}$, and on distinct planes, the contraction by
$S$ will be equal to $N$, as required.  In the
indicated free product, $\{a,b\}$ and $\{c,d\}$ are simply
``in general position'' on such planes.

We prove next that free product respects matroid duality and is
associative.  First, recall that for any matroid $M(S)$, the rank
function of the dual matroid $M\dual$ satisfies $\rkmd B = |B|-\rho
(M)+\rkm A$, or equivalently, $\rlm A=\nlmd B$, whenever $A+B=S$.
\begin{pro}[\cite{crsc:fpm}]\label{pro:dual}
  For all matroids $M$ and $N$, $(M\frp N)\dual = N\dual\frp M\dual$.
\end{pro}
\begin{proof}
Suppose that $M=M(S)$, $N=N(T)$, and $A+B=S+T$, so that $A$ is
a basis for $M\frp N$ if and only if $B$ is a basis for $(M\frp N)\dual$.
Then $A$ is a basis for $M\frp N$ if and only if \as\ is independent
in $M$, \at\ spans $N$ and $\rlm\as=\nln\at$, which is true if and only
if \bs\ spans $M\dual$, \bt\ is independent in $N\dual$, and 
$\nlmd\bs=\lambda_{\scriptscriptstyle N^*}(\bt)$, that is, if
and only if $B$ is a basis for $N\dual\frp M\dual$.
\end{proof}

\begin{pro}
Free product is an associative operation.
\end{pro}
\begin{proof}
Suppose that $M=M(S)$, $N=N(T)$ and $P=P(U)$.  A set $A\subseteq S+T+U$
is independent in $(M\frp N)\frp P$ if and only if $A\sub{S+T}$ is 
independent in $M\frp N$ and 
$\lambda_{\scriptscriptstyle M\frp N}(A\sub{S+T})\geq\nlp\au$.  
Since $A\sub{S+T}$ is independent in $M\frp N$,
we have
\begin{align*}
  \lambda\sub{M\frp N}(A\sub{S+T}) &\= \rho (M\frp N)-|A\sub{S+T}|\\
&\= \rho (M) +\rho (N) - |\as|-|\at|\\
&\= \rlm\as +\rho (N)-|\at|.
\end{align*}
Hence $A$ is independent in $(M\frp N)\frp P$ if and only if
$\as$ is independent in $M$,
$\nln\at\leq\rlm\as$ and $\nlp\au\leq\rlm\as +\rho (N)-|\at|.$
Adding $\nln\at$ to both sides of the last inequality, we may
express these three conditions as
$$
\nlm\as\leq 0,\qquad\nln\at\leq\rlm\as\sands\nln\at +\nlp\au
\leq\rlm\as +\rln\at.
$$
On the other hand,  $A$ is independent in $M\frp (N\frp P)$ if and only if
$\nlm\as\leq 0$ and $\nu\sub{N\frp P}(A\sub{T+U})\leq
\rlm\as$.  By Equation \ref{eq:null}, the latter inequality
may be written as
$$
\nln\at +\nlp\au\leq\rlm\as +\min\{\rln\at,\nlp\au\},
$$
which holds if and only if
$\nln\at\leq\rlm\as$ and $\nln\at +\nlp\au\leq\rlm\as +\rln\at$.
Hence $A$ is independent in $M\frp (N\frp P)$ if and only if
it is independent in $(M\frp N)\frp P$.
\end{proof}

The definitions and properties stated above have natural
analogs for iterated free products.  
\begin{pro}\label{pro:nindep}
  If $L(S)=M_1(S_1)\frp\cdots\frp M_k(S_k)$, then 
$A\subseteq S$ is independent in $L$ if and only if
\begin{equation}
  \label{eq:nindep}
\sum_{i=1}^{j-1}\lambda\sub{M_i}(A\sub{S_i})\geq
\sum_{i=1}^j\nu\sub{M_i}(A\sub{S_i}),
\end{equation}
for all $j$ such that $1\leq j\leq k$. 
\end{pro}
\begin{proof}  
We use induction on $k$.
When $k=1$, the sum on the left-hand side of the inequality
is empty and thus zero; so the result holds.  Suppose the
result holds for $L'=M_1(S_1)\frp\cdots\frp 
M_{k-1}(S_{k-1})$.  Then $A$ is independent in $L=L'\frp M_k$
if and only if $A\sub{S_k}'=A\sub{S_1}+\cdots+A\sub{S_{k-1}}$ is independent
in $L'$ and $\nu\sub{M_k}(A\sub{S_k})\leq\lambda\sub{L'}(A\sub{S_k}')$,
that is, if and only if Inequality \ref{eq:nindep} holds
for $1\leq j\leq k-1$ and, since $A\sub{S_k}'$ is independent in $L'$, 
$$
\nu\sub{M_k}(A\sub{S_k})\leq\rho (L')-|A\sub{S_k}'|
\= \sum_{i=1}^{k-1}\rho (M_i)-|A\sub{S_i}|.
$$
But $\rho(M_i)-|A\sub{S_i}|=
\lambda\sub{M_i}(A\sub{S_i})-\nu\sub{M_i}(A\sub{S_i})$, for
all $i$; hence the above inequality is equivalent to
Inequality \ref{eq:nindep}, for $j=k$.
\end{proof}
We will need the following generalization of Proposition
\ref{pro:univ} in Section \ref{sec:algebra}.
\begin{pro}\label{pro:guniv}
  Suppose that $L=L(U)$ and $\emptyset=T_0\subset\cdots\subset
T_k=U$ is a chain of subsets of $U$, for some $k\geq 0$, and  
let $L_i$ denote the minor $L(T_{i-1},T_i)$, for $1\leq i\leq k$.
The identity map on $U$ is a weak map $L_1\frp\cdots\frp L_k\rta
L$.
\end{pro}
\begin{proof}
  Let $S_i=T_i\backslash T_{i-1}$, for $1\leq i\leq k$, so
that $L_i=L_i(S_i)$, for all $i$.  By Lemma \ref{lem:conineq} and
induction on $k$, it follows that the inequalities \eqref{eq:nindep} 
hold for all independent sets $A$ in $L$.  Hence, by 
Proposition \ref{pro:nindep}, any independent set in $L$ is also independent
in $L_1\frp\cdots\frp L_k$, that is, the identity map on $U$ is a weak
map $L_1\frp\cdots\frp L_k\rta
L$.
\end{proof}

One-element matroids (isthmuses and loops) play a special
role in the study of free products.  

\begin{exa}\label{exa:unif}
Recall that, if $\{a\}$ is any  singleton, then $I(a)$ and
$Z(a)$ denote the matroids on $\{a\}$ consisting, respectively,
of a single point and a single loop.
For any set $S=\{s_1,\dots,s_n\}$, and $k\leq n$, the
free product $I(s_1)\frp\cdots\frp I(s_k)\frp Z(s_{k+1})\frp
\cdots\frp Z(s_n)$ is the uniform matroid $U_{k,n}(S)$.
\end{exa}

For any matroid $M$, we write \loo M\ and \ist M, respectively, for the 
sets of loops and isthmuses of $M$.  
\begin{pro}\label{pro:ptlpcon}
  For all matroids $M$ and $N$, $\loo M\subseteq\loo{M\frp N}$,
with $\loo M = \loo{M\frp N}$, whenever $\rho (M)>0$.  Dually,
$\ist N\subseteq\ist{M\frp N}$, with equality whenever $\nu (N)>0$.
\end{pro}
\begin{proof}
  If $x$ is a loop of $M$, then $x$ belongs to no independent
set of $M\frp N$; hence $x$ is a loop of $M\frp N$, and so
$\loo M\subseteq\loo{M\frp N}$.  On the other hand, suppose that
$\rho (M)>0$, and that $N=N(T)$ and $x\in T$. It follows from 
Proposition \ref{pro:rank} that 
$\rho\sub{M\frp N}(x)=\rkn x +\min\{\rho (M),
\nln x\}=1$, so $x$ is not a loop in $M\frp N$, and hence
$\loo{M\frp N}=\loo M$.
The dual statements follow directly from Proposition \ref{pro:dual}.
\end{proof}

\begin{cor}\label{cor:ptlp}
  If $\rho (M)=0$ or $\nu (N)=0$, then $M\frp N=M\oplus N$.
\end{cor}

\begin{exa}
 For any matroid $M$, the matroids $M\frp I$
and $Z\frp M$ consist of $M$ with, respectively, an isthmus and
a loop adjoined, while $M\frp Z$ and $I\frp M$ are
respectively the free one-point extension and coextension of $M$
(see \cite{ox:mt}).
\end{exa}

\begin{exa}\label{exa:freedom}
Because adjoining an isthmus and taking a single-point free extension
of a matroid correspond to free multiplication on the right by $I$ and
$Z$, respectively, it follows that the class of matroids introduced in
\cite{cr:see}, now variously known as {\it generalized Catalan matroids} 
(\cite{bodeno:lpm}), {\it shifted matroids} (\cite{ar:tcm}) and
{\it freedom matroids} (\cite{crsc:fsa}), is the class generated by the 
single-element matroids under free product.
\end{exa}

A {\it representation} of a matroid $M(S)$ over a field $F$ is a
matrix $P$ having entries in $F$ and rows labeled by the elements of
$S$, such that for all $A\subseteq S$, the
submatrix $P\sub A$ of $P$, consisting of those rows of $P$ whose
labels belong to $A$, has rank $\rkm A$.  We can, and shall, always
assume that the number of columns in a representation of $M$ is equal
to the rank of $M$.  A matroid $M$ is called $F$-{\it representable}
if there exists a representation of $M$ over $F$.

\begin{pro}\label{pro:repprod}
  If the matroids $M(S)$ and $N(T)$ are $F$-representable,
and the field $F$ is large enough, then the free product
$M\frp N$ is $F$-representable. 
\end{pro}
\begin{proof}
Suppose that $P$ and $Q$ are representations for $M$ and $N$,
respectively. Using the fact that the field $F$ has enough elements,
we can construct a $|T|\times\rho (M)$ matrix $Z$, with rows labelled
(arbitrarily) by $T$, having the following 
property:  given any $A\subseteq S$ which is independent in $M$,
and any $B\subseteq T$ of size $\rlm A =\rho(M)-|A|$, the matrix
$$
\left[
\begin{array}{c}
\raisebox{0ex}[2.35ex][1ex]{\!\!$P\sub A$\!\!}\\ \hline
\raisebox{0ex}[2.35ex][1ex]{\!\!$Z\sub B$\!\!}
\end{array}
\right]
$$
is nonsingular. We show that the matrix
$$
R\= \left[
\begin{array}{c|c}
\raisebox{0ex}[2.35ex][1ex]{$P$}
 & 0\\ \hline
\raisebox{0ex}[2.35ex][1ex]{$Z$} & Q
\end{array}
\right]
$$
is a representation for the free product $M\frp N$.
Suppose that $A\subseteq S+T$, and let $B\subseteq\at$ be a 
basis for \at\ in $N$.  Since $B$ is independent
in $N$, the matrix $Q\sub B$ has independent rows, and
hence the matrix $R\sub A$ has independent rows if and only
if the matrix
$$
\left[
\begin{array}{c}
\raisebox{0ex}[2.35ex][1ex]{\!\!$P\sub\as$\!\!}\\ \hline
\raisebox{0ex}[2.35ex][1ex]{\!\!$Z\sub{\at\backslash B}$\!\!}
\end{array}
\right]
$$
has independent rows.  Since $|\at\backslash B|=\nln\at$,
it follows from the construction of $Z$ that
this latter matrix has independent rows if and only if
$\as$ independent in $M$ and $\rlm\as\geq\nln\at$, that is,
if and only if $A$ is independent in $M\frp N$.
\end{proof}

Suppose that $\ca=\{A_i\co i\in I\}$ is an indexed family of subsets
of a set $S$ (with repetitions allowed). A set $A\subseteq S$ is a
{\it partial transversal} of \ca\ if there exists an injective map
$f\co A\rta I$ such that $a\in A_{f(a)}$, for all $a\in A$.  The set
of partial transversals of \ca\ is the collection of independent sets
of a matroid, called a {\it transversal matroid} on $S$, and denoted
by $M(S, \ca)$.  The family \ca\ is a {\it presentation} of 
$M(S,\ca)$.  Any transversal matroid $M$ has a presentation
with number of sets equal to the rank of $M$ (see \cite{we:mt}, page
244).

\begin{pro}\label{pro:trans}
The free product of transversal matroids is a transversal matroid.
\end{pro}
\begin{proof}
Suppose that $M=M(S,\ca)$ and $N=M(T,\cb)$ are transversal matroids with
respective presentations $\ca = \{A_i\co i\in I\}$ and $\{ B_j\co j\in J\}$, 
where $|I|=\rho (M)$.  
For all $k\in I+J$, define $U_k\subseteq S+T$ by
$$
U_k \=
\begin{cases}
  A_k +T & \text{if $k\in I$,}\\
  B_k & \text{if $k\in J$.}
\end{cases}
$$
We show that the free product $M\frp N$ is equal to
the transversal matroid on $S+T$ having presentation
$\cu =\{U_k\co k\in I+J\}$.  
Given $A\subseteq S+T$, let $B\subseteq\at$ be a basis for 
\at\ in $N$.
The set $A$ is independent in $M(S+T,\, \cu)$
if and only if there exists injective $f\co A\backslash B\rta I$
such that $a\in U_{f(a)}$ for all $a\in A\backslash B$, which is
the case if and only if \as\ is independent
in $M$ and $|\at\backslash B|\leq |I|-|\as|$.  Since
$|\at\backslash B|=\nln\at$ and $\rlm\as = |I|-|\as|$,
for \as\ independent in $M$, it follows that such $f$ exists
if and only if $A$ is independent in $M\frp N$.
\end{proof}

\section{Minors of free products}\label{sec:minors}
The minors of a free product of matroids are perhaps most
simply described in terms of the matroid truncation operator 
and its dual, the Higgs lift operator (see \cite{hi:smg}).
The {\it truncation\/} of a matroid $M(S)$
is the matroid $\trun{}M$ whose independent sets are those
independent sets $A$ of $M$ satisfying $|A|\leq\max\{0,\rho(M)-1\}$,
and the {\it Higgs lift\/}, or simply {\it lift}, of $M$
is the matroid $\lift{}M$ whose family of independent sets
is $\{ A\subseteq S\coo\nlm A\leq 1\}$. Denoting by
$\trun iM$ and $\lift iM$, respectively, the $i$-fold truncation
and lift of $M(S)$, it follows that
$\trun i M$ has rank equal to $\max\{ 0,\rho (M)-i\}$, and
$$
\rho\sub{\trun iM}(A)\=\min\{\rkm A,\,\rho (\trun iM)\}\sands
\lambda\sub{\trun iM}(A)\=\min\{0,\,\rlm A-i\},
$$
for all $A\subseteq S$. The rank of $\lift i M$
is $\min\{ |S|,\rho (M)+i\}$, and
$$
\rho\sub{\lift iM}(A)=\min\{|A|,\,\rkm A+i\}\sands\
\nu\sub{\lift iM}(A)=\max\{0,\,\nlm A -i\}
$$
for all $A\subseteq S$.  The truncation and lift operators
are dual to each other, so that $(\trun i M)\dual
=\lift i (M\dual)$, for all matroids $M$ and $i\geq 0$. Truncation
commutes with contraction and lift commutes
with restriction, so for any matroid $M(S)$ and $i\geq 0$,
$$
(\trun i M)/U\=\trun i (M/U)\spandsp
(\lift i M)|U\= \lift i (M|U),
$$
for all $U\subseteq S$.  We thus shall write expressions such as
these without parentheses.  The precise manner in which lift and
truncation fail to commute with contraction and restriction,
respectively, is described by the following proposition.
\begin{pro}\label{pro:ltnc}
For any matroid $M(S)$ and $U\subseteq S$
$$
  \trun i (M|U)\= (\trun{i+j}M)|U\spandsp
\lift i (M/U)\= (\lift{i+k}M)/U,
$$
for all $i\geq 0$, where $j=\rlm U$ and $k=\nlm U$.
\end{pro}
\begin{proof}
  The rank-lack of $A\subseteq U$ in $M|U$ is given by
$\lambda\sub{M|U}(A)=\rlm A-\rlm U=\rlm A-j$, and so
$\lambda\sub{\trun i{(M|U)}}(A) =
\min\{0,\,\lambda\sub{M|U}(A)-i\}=\min\{0,\rlm A-j-i\}$.
On the other hand, 
\begin{align*}
\lambda\sub{(\trun{i+j}M)|U}(A)&\= \lambda\sub{\trun{i+j}M}(A)-
\lambda\sub{\trun{i+j}M}(U)\\
&\= \min\{ 0,\,\rlm A - i-j\} -\min\{ 0,\,\rlm U - i - j\},
\end{align*}
which is equal to $\min\{0,\rlm A - i -j\}$, since $\rlm U = j$.
The matroids $\trun i (M|U)$ and $(\trun{i+j}M)|U$ thus have identical
 rank-lack functions, and are therefore equal.  
The second equality follows
from duality, using the fact that $\rlm U=\nu\sub{M\dual}(S\backslash U)$,
for all $U\subseteq S$.
\end{proof}

In keeping with the notational tradition of performing unary
operations before binary operations, in order to avoid a proliferation
of parentheses, we adopt the convention that
all truncations, lifts, deletions and contractions that may appear in
a given expression for a matroid are to be performed before any free
products and/or direct sums that appear.

\begin{pro}\label{pro:delcon}
  If $P=M(S)\frp N(T)$ and $U\subseteq S+T$, then
$$
P|U\= M|\us\frp\lift i N|\ut\spandsp
P/U\= \trun j M/\us \frp N/\ut,
$$
where $i=\rlm\us$ and $j=\nln\ut$.
\end{pro}
\begin{proof}
  A set $A\subseteq U$ is independent in $P|U$ 
if and only if \as\ is independent in $M$ and $\rlm\as\geq\nln\at$.
Using the fact that $\rlm\as = \lambda\sub{M|\us}(\as)+\rlm\us$
and that $\nln\at =\nu\sub{N|\ut}(\at)$, we
thus have $A$  independent in $P|U$ if and only if
\as\ is independent in $M|\us$ and
$\lambda\sub{M|\us}(\as)\geq\nu\sub{N|\ut}(\at)-i$.
But $\max\{ 0, \nu\sub{N|\ut}(\at)-i\}
=\nu\sub{\lift i N|\ut}(\at)$, and so
$A$ is independent in $P|U$ if and only if $\as$ is independent in
$M|\us$ and $\lambda\sub{M|\us}(\as)\geq\nu\sub{\lift i N|\ut}(\at)$,
that is, if and only if $A$ is independent in 
$M|\us\frp\lift i N|\ut$.

The second equality follows from the first by duality, that is,
by Proposition \ref{pro:dual}, the duality between deletion
and contraction, the duality between lift and truncation
and the fact that $\lambda\sub{N\dual}(T\backslash\ut)=\nln\ut$.
\end{proof}
\begin{thm}\label{thm:minors}
  If $P=M(S)\frp N(T)$ and $U\subseteq V\subseteq S+T$, then
$$
P(U,V)\= (\trun j M)(\us,\vs )\frp (\lift i N)(\ut,\vt),
$$
where $j=\nln\ut$ and $i=\rlm\vs$.
\end{thm}
\begin{proof}
  By Proposition \ref{pro:delcon}, we have
$P|V = M|\vs\frp (\lift i N|\vt)$, where $i=\rlm\vs$, and
thus, by the same proposition,
\begin{align*}
P(U,V) \= (P|V)/U &\= (\trun k (M|\vs))/\us \frp (\lift i N|\vt)/\ut\\
&\=(\trun k (M|\vs))/\us \frp (\lift i N)(\ut,\vt),
\end{align*}
where $k=\nu\sub{\lift i N|\vt}(\ut)=\max\{0, \nln\ut-i\}=
\max\{0,j-i\}$. If $j\geq i$, then by Proposition \ref{pro:ltnc}, 
\begin{align*}
(\trun k (M|\vs))/\us& \= ((\trun{k+i}M)|\vs)/\us\\
&\= (\trun jM)(\us,\vs),
\end{align*}
and we thus obtain the desired expression for $P(U,V)$.  On the other
hand, if $j<i=\rlm\vs$, then $(\trun jM)|\vs= M|\vs$, and
$k=0$, and thus 
\begin{align*}
(\trun k (M|\vs))/\us& \= (M|\vs)/\us\\
&\= ((\trun j M)|\vs)/\us\\
&\= (\trun jM)(\us,\vs),
\end{align*}
and again we obtain the desired expression for $P(U,V)$.
\end{proof}

As a special case of Theorem \ref{thm:minors}, we have that
the minors of $P=M(S)\frp N(T)$  supported on
the sets $S$ and $T$ are obtained by successive truncations of $M$ and  
Higgs lifts of $N$, respectively; that is,
for all $A\subseteq S$ and $B\subseteq T$,
$$
P(A,A\cup T)\= \lift i N \spandsp P(B,B\cup S)\=\trun j M,
$$
where $i=\rlm A$ and $j=\nln B$.
This is to be compared to the
direct sum, where these minors are simply isomorphic to $M$ and $N$.

The following proposition describes how the lift and truncation
operators interact with free product.

\begin{pro}\label{pro:ltfrp}
  For all matroids $M$ and $N$,
The truncation and lift of the free product $M\frp N$ are given by
$$
\trun{}{(M\frp N)}\=
\begin{cases}
  M\frp \trun{}N & \text{if $\rho (N)>0$,}\\
 \trun{}M\frp N & \text{if $\rho (N)=0$,}
\end{cases}
$$
and
$$
\lift{}{(M\frp N)}\=
\begin{cases}
  \lift{}M\frp N & \text{if $\nu (M)>0$,}\\
 M\frp\lift{}N & \text{if $\nu (M)=0$,}
\end{cases}
$$
for all matroids $M$ and $N$.
\end{pro}
\begin{proof}
  If $\rho (M)=0$ then, by Corollary \ref{cor:ptlp}, we have
 $M\frp N = M\oplus N$, and so $\trun{}(M\frp N)=
\trun{}(M\oplus N) =\trun{}M\oplus\trun{}N= M\oplus\trun{}N= M\frp
\trun{}N$.  We therefore assume that $\rho (M)$ is nonzero.

Suppose that $M=M(S)$ and $N=N(T)$.  Observe that if a set
$A\subseteq S+T$ is independent in any of the matroids
$\trun{}(M\frp N)$, $M\frp\trun{}N$ and $\trun{}M\frp N$, then
\as\ is necessarily independent in $M$. Hence, for the remainder of
the proof, we assume that $A$ is some subset of $S+T$ such that
\as\ is independent in $M$.

We first consider the case in which $\rho (N)=0$.  The set $A$ is
independent in $M\frp N$ if and only if $\rlm\as\geq\nln\at$, which
is the case if and only if $|A|\leq\rho(M)$, since
$\rlm\as= \rho (M)-|\as|$ and $\nln\at=|\at|$.  It follows that
$A$ is independent in $\trun{}(M\frp N)$ if and only if $|A|\leq
\rho(M)-1$.  

Now $A$ is independent in $\trun{}M\frp N$ if and only if
$\rkm\as=|\as|\leq\rho (M)-1$ and $\lambda\sub{\trun{}M}(\as)\geq
\nln\at$. Furthermore
\begin{align*}
  \lambda\sub{\trun{}M}(\as)& \=\max\{\rlm\as-1,\,0\}\\
&\= \max\{\rho(M)-|\as|-1,\,0\},
\end{align*}
which is equal to $\rho(M)-|\as|-1$, since $|\as|\leq\rho(M)-1$.
Therefore $A$ is independent in $\trun{}M\frp N$ if and only if
$\rho(M)-|\as|-1\geq\nln\at= |\at|$, that is, if and only if 
$|A|\leq\rho (M)-1$, and hence $\trun{}(M\frp N)=\trun{}M\frp N$.

Now suppose that $\rho (N)>0$.  If $\rkn\at<\rho(N)$ then,
by Proposition \ref{pro:bases}, the set
$A$ doesn't span $M\frp N$, and so $A$ is independent in $\trun{}(M\frp
N)$ if and only if $A$ is independent in $M\frp N$.  But since
\at\ doesn't span $N$, and thus $\nu\sub{\trun{}N}(\at)=\nln\at$,
it follows that $A$ is independent in $M\frp N$ if and only if it is
also independent in $M\frp\trun{}N$.  
If $\rkn\at=\rho(N)$ then, by Proposition \ref{pro:bases}, we have that
$A$ is independent in $\trun{}(M\frp N)$ if and only if
$\rlm\as>\nln\at$.  
But $A$ is independent in $M\frp\trun{}N$
if and only if $\rlm\as\geq\nu\sub{\trun{}N}(\at)=\nln\at+1$;
hence $\trun{}(M\frp N)=M\frp\trun{}N$.
The corresponding result for  $\lift{}(M\frp N)$ follows by duality.
\end{proof}

\noindent
It follows from Proposition \ref{pro:ltfrp} that, for all
matroids $M$ and $N$, and $i\geq 0$,
\begin{equation}\label{eq:ltfrp}
\trun i{(M\frp N)}\= \trun jM\frp\trun{i-j}N \sands
\lift i{(M\frp N)}\= \lift{i-k}M\frp\lift kN, 
\end{equation}
where $j=\max\{i-\rho(N),0\}$ and $k=\max\{i-\nu (M),0\}$.

\section{Irreducible matroids and unique factorization}\label{sec:uft}

A crucial tool for the study of factorization of matroids with respect
to free product is the notion of {\it cyclic flat} of a matroid.
Recall that a cyclic flat of $M$ is a flat $A$ which is
equal to a union of circuits of $M$.  Alternatively, a flat $A$ is
cyclic if and only if the restriction $M|A$ is isthmusless.  Observe
that in particular, any closure of a circuit in a matroid is a cyclic
flat.  We begin with the following characterization of the cyclic
flats in a free product of matroids.
\begin{pro}\label{pro:cffrp}
  A subset $A\neq S$ of $S+T$ is a cyclic flat of 
$L=M(S)\frp N(T)$ if and only if either $A\subseteq S$ and
$A$ is a cyclic flat of $M$, or $A=S\cup B$, where $B$ is
a (nonempty) cyclic flat of $N$.  The set $S$ is a cyclic flat
of $L$ if and only if $M$ is isthmusless and $N$ is loopless.
\end{pro}
\begin{proof}
  Suppose that $A\subseteq S+T$ satisfies $\rlm\as>\nln\at$ and
$A\neq S$. According to Corollary \ref{cor:flats}, $A$ is a flat of
$L$ if and only if \as\ is a flat of $M$, in which case
any element of \at\ is an isthmus of $L|A$.  
Hence $A$ is a cyclic flat of $L$ if and only if 
$\at=\emptyset$ and $A=\as$ is a cyclic flat of $M$.

Now suppose that $A\neq S$ and $\rlm\as\leq\nln\at$. Then by 
Corollary \ref{cor:flats},
$A$ is a flat of $L$ if and only if $\as=S$ and \at\ is a nonempty
flat of $N$.  Given such a flat $A$, we have
$\rkl A = \rkm\as +\rkn\at +\min\{\rlm\as,\nln\at\}=\rho (M)+\rkn\at$;
hence if $A$ is cyclic then \at\ must be a cyclic flat of $N$.  
On the other hand,
if \at\ is cyclic in $N$, then $\rkl{A\backslash a}=\rkl A$,
for all $a\in\at$, and since $\nln\at>0$ and $\rlm\as=\rlm S=0$, it
follows that $\rkl{A\backslash a}=\rkl A$ for all $a\in\as$ as well.
Hence $A$ is cyclic.

Since $\rlm S = 0$, it follows from Corollary \ref{cor:flats}  that
$S$ is a flat of $L$ if and only if $N$ is loopless, in which case
the flat $S$ is cyclic if and only if $M=L|S$ is isthmusless.  

\end{proof}

\begin{dfn}{\rm
A set $A\subseteq S$ is a {\it free separator\/} of a
matroid $M(S)$ if every cyclic flat of $M$ is comparable to $A$ by
inclusion.
}    
\end{dfn}
Note that the empty set and the entire set $S$ are free
separators of any matroid $M(S)$; any other free separator is said to
be {\it nontrivial}.

\begin{thm}\label{thm:freesep}
  For any matroid $L(S+T)$, the following are equivalent:
  \begin{enumerate}[(i)]
\item $L(S+T)\= L|S\frp L/S$.
\item $S$ is a free separator of $L$.
 \end{enumerate}
\end{thm}
\begin{proof}
The implication $(i)\Rightarrow (ii)$ is immediate from
Proposition \ref{pro:cffrp}.
  Conversely, suppose that $S$ is a free separator of $L$, and
let $M=L|S$ and $N=L/S$.  We first show that every circuit of $L$
is also a circuit of the free product $M(S)\frp N(T)$.  Let $C$ be a circuit
of $L$.  If $C\subseteq S$, then $C$ is a circuit of $M$, and therefore
a circuit of $M\frp N$. Suppose that $C\not\subseteq S$. Since
$C$ is a circuit, $\rkl{C\backslash a} = \rkl C$
and thus, by the semimodularity of the rank function,
$\rkl{(S\cup C)\backslash a}=\rkl{S\cup C}$,  for all $a\in C$.
Hence, for all $a\in\ct$, we have $\rkn{\ct}=\rkl{S\cup C}-\rkl S=
\rkl{S\cup C\backslash a}-\rkl S=\rkn{\ct\backslash a}$,
and so $N|\ct$ is isthmus free.
Since the closure of a circuit is a cyclic flat, $S$ is a free
separator, and $C\not\subseteq S$,
we have $S\subseteq\cll C$. It follows that $\rkl{S\cup C}=\rkl C =
|C|-1$, and so $\nll{S\cup C}= |S|-|\cs|+1$.  Therefore
\begin{align*}
\nln\ct & \= \nll{S\cup C} - \nll S\\
& \= |S|-|\cs|+1 -(|S|-\rkl S)\\
&\= \rho (M) - |\cs| +1,
\end{align*}
which is equal to $\rlm\cs +1$, since $\cs$ is independent
in $L$ (and thus also in $M$).  By Proposition \ref{pro:circuit},
it follows that $C$ is a circuit in $M\frp N$. 

We have thus shown that every circuit in $L$ is also a circuit in
$L|S\frp L/S$, in other words, the identity map on $S+T$ is a weak
map $L\rta L|S\frp L/S$.  By \ref{pro:univ}, the identity map on 
$S+T$ is also a weak map $L|S\frp L/S\rta L$; hence $L=L|S\frp L/S$. 
\end{proof}

We refer to a nonempty matroid $M$ as {\it irreducible} if any
factorization of $M$ as a free product of matroids contains $M$ as a
factor.  By convention, the empty matroid is not irreducible.  
The following restatement of Theorem \ref{thm:freesep} characterizes
irreducible matroids.
\begin{thm}\label{thm:irred}
  For any nonempty matroid $M(S)$, the following are equivalent:
  \begin{enumerate}[(i)]
  \item $M$ is irreducible with respect to free product.
\item $M$ has no nontrivial free separator.
  \end{enumerate}
\end{thm}
\begin{cor}\label{cor:ds}
   If $M$ is loopless, isthmusless and disconnected, then $M$
is irreducible.  
\end{cor}
\begin{proof}
  Suppose that $M(S)$ is loopless, isthmusless and disconnected,
and write $M(S)$ as the direct sum $P(U)\oplus Q(V)$, with $U$ and
$V$ nonempty. Let $A$ be a nonempty 
proper subset of $S$.  Assume, without loss of generality,
that \au\ and $V\backslash A$ are nonempty, and let
$a\in V\backslash A$.  Since $Q$ is loop and isthmus
free, $a$ is contained in some circuit $C$ of $Q$.  Now $C$
is also a circuit of $M$ and $a\in \clm C=\clq C\subseteq V$;
hence $\clm C$ neither contains nor is contained in $A$, and
so $A$ is not a free separator of $M$.
\end{proof}

\begin{cor}\label{cor:unif}
  If $L=M(S)\frp N(T) = P(T)\frp Q(S)$, where $S$ and $T$ are
nonempty, then $L$ is a uniform matroid.
\end{cor}
\begin{proof}
  Let $C$ be a circuit of $L$.  
By Theorem \ref{thm:freesep}, both $S$ and $T$ are free
separators of $L$ and hence $\cll C$ is comparable to both
$S$ and $T$ by inclusion. Since $S$ and $T$ are disjoint and nonempty,
the only possibility is that $S$ and $T$ are both contained
in $\cll C$.  Every circuit of $L$ is thus a spanning set for $L$,
and therefore $L$ is uniform.
\end{proof}

We remark 
that it follows from Proposition \ref{pro:dual} that a matroid $M$ is
irreducible if and only if the dual matroid $M\dual$ is irreducible.
\begin{cor}\label{cor:sd}
  If $M$ is identically self-dual, then $M$ is either uniform
or irreducible.
\end{cor}
\begin{proof}
Suppose that $M$ is identically self-dual and factors
as $P(U)\frp Q(V)$, with $U$ and $V$ nonempty.  Using Proposition
\ref{pro:dual}, we have $P(U)\frp Q(V) = M= M\dual = Q\dual (V)\frp
P\dual (U)$, and hence it follows from Corollary \ref{cor:unif} that
$M$ is uniform.  
\end{proof}
\begin{exa}
  Suppose that $S=\{a,b,c,d\}$ and let $M(S)$ be the matroid
in which $ab$ is a double point, collinear with $c$ and $d$.
Then $M$ is self-dual, not uniform, and factors with respect to
free product as
$I(a)\frp Z(b)\frp I(c)\frp Z(d)$.
\end{exa}

For any matroid $M(S)$, we denote by $\dlat (M)$ the complete
sublattice of the Boolean algebra $2^S$ generated by all cyclic flats
of $M$. Note that $\dlat (M)$ is a distributive lattice, and contains
in particular the empty union and empty intersection of cyclic flats
of $M$, which are equal to $\emptyset$ and $S$, respectively.

\begin{pro}\label{pro:unif}
  A nonempty matroid $M(S)$ is uniform if and only if $|\dlat (M)|=2$,
that is, if and only if $\dlat (M)=\{\emptyset, S\}$.
\end{pro}
\begin{proof}
  Uniform matroids are characterized by the fact that
all of their circuits are spanning. Hence $M(S)$ is
uniform if and only if it has no cyclic flat that is
both nonempty and not equal to $S$.  For nonempty matroids,
this is the case if and only if $\dlat (M)=\{\emptyset, S\}$.
\end{proof}
\begin{dfn}{\rm
An element $x$ of a partially ordered set $P$ is a
{\it pinchpoint} if the set $\{x\}$ is a crosscut of $P$,
that is, if all elements of $P$ are comparable to $x$.
A pinchpoint of $P$ is {\it nontrivial} if it is neither 
minimal nor maximal in $P$.
}\end{dfn}

A uniform matroid is irreducible with respect to free product if and
only if its underlying set is a singleton (see Example
\ref{exa:unif}).  Irreducibility of nonuniform
matroids is characterized in the following theorem.  

\begin{thm}\label{thm:irred2}
  For any nonuniform matroid $M(S)$, the following are equivalent:
  \begin{enumerate}[(i)]
  \item $M$ is irreducible with respect to free product.
\item The lattice $\dlat (M)$ contains no nontrivial pinchpoint.
  \end{enumerate}
\end{thm}
\begin{proof}
If $A\in\dlat (M)$ is a nontrivial pinchpoint then $A\subseteq S$ is itself
a nontrivial free separator, and hence $M$ is not irreducible by 
Theorem \ref{thm:irred}.  Conversely, suppose that
$M(S)$ is nonuniform and has a nontrivial free separator
$A\subseteq S$.  Since $M$ is nonuniform it has a cyclic flat
$B$ which is neither empty nor equal to $S$.  
If $A\subseteq B$, then the intersection
of all cyclic flats of $M$ containing $A$ is a nontrivial pinchpoint of 
$\dlat (M)$.  If $B\subseteq A$, then the union of all
cyclic flats which are contained in $A$ is a nontrivial pinchpoint.
\end{proof}

For any matroid $M(S)$ we denote by $\fs (M)$ the set of all free
separators of $M$, ordered by inclusion.  We shall see presently that
$\fs (M)$ is a lattice (in fact distributive).
For all $A\subseteq B\subseteq S$, we denote by $[A,B]$ the
subinterval $\{ U\subseteq S\co A\subseteq U\subseteq B\}$ of the
Boolean algebra $2^S$.  If $A$ and $B$ are free separators of $M(S)$,
then we write $[A,B]\sub\cf$ for the subinterval $[A,B]\cap\fs (M)$ of
$\fs (M)$.  In the following lemma we show that an interval in the
lattice of free separators of a matroid is isomorphic, under the
obvious map, to the lattice of free separators of the corresponding
minor of the matroid.

\begin{lem}\label{lem:fspint}
  For all free separators $A\subseteq B$ of a matroid $M(S)$,
the map from the interval $[A,B]\sub\cf$ in $\fs (M)$ to the
lattice $\fs (M(A,B))$ given by $U\mapsto U\backslash A$
is a bijection (and thus a lattice isomorphism).
\end{lem}
\begin{proof}
  If $A\subseteq U\subseteq B$ are free separators of $M(S)$,  
then it follows from Theorems \ref{thm:minors} and \ref{thm:freesep} that 
$M(A,B)=M(A,U)\frp M(U,B)$, and so $U\backslash A$
is a free separator of $M(A,B)$.  On the other hand, 
if $A\subseteq B$ are free separators of $M$, then $M$ 
factors as $M= M|A\frp M(A,B)\frp M/B$, and if $V\subseteq
B\backslash A$ is a free separator of $M(A,B)$, we have the factorization
$M(A,B)=M(A,B)|V\frp M(A,B)/V = M(A, A\cup V)\frp M(A\cup V,B)$.
Hence, by associativity of free product, $A\cup V$ is
a free separator of $M$.
\end{proof}

If $U_0\subset\cdots\subset U_k$ is a chain in $\fs (M)$, with
$U_0=\emptyset$ and $U_k=S$, then by Lemma \ref{lem:fspint}, 
we have the factorization
$
M(S)\= M(U_0,U_1)\frp\cdots\frp M(U_{k-1},U_k)
$
of $M$ into a free product of nonempty matroids.  On the
other hand, given any factorization $M(S)=M_1(S_1)\frp\cdots\frp M_k(S_k)$,
with all $M_i$ nonempty, the sets $U_i=S_0\cup\cdots\cup S_i$, for
$1\leq i\leq k$, comprise a chain from $\emptyset$ to $S$ in $\fs (M)$.  
Hence the factorizations of $M(S)$ into free products of nonempty
matroids are in one-to-one correspondence with chains from $\emptyset$
to $S$ in the lattice $\fs (M)$.

\begin{lem}\label{lem:fspuni}
  A matroid $M(S)$ is uniform if and only if $\fs (M)$ is
equal to the Boolean algebra $2^S$.
\end{lem}
\begin{proof}
  If $M(S)$ is uniform then the only possible cyclic flats
of $M$ are $\emptyset$ and $S$, and so every subset of $S$ is
a free separator of $M$.  Conversely, if every subset of $S$ is
a free separator of $M$, then the only possible cyclic flats
of $M$ are $\emptyset$ and $S$, and thus $M$ must be uniform.
\end{proof}

\begin{dfn}{\rm
  The {\it primary flag} $\calt\sub M$ of a matroid $M$ is the
chain $T_0\subset\cdots\subset T_k$ consisting of all pinchpoints
in the lattice $\dlat (M)$.
}
\end{dfn}
Note that the sets belonging to the primary flag of a matroid are, in
particular, free separators, and thus the primary flag
of $M$ is a chain from $\emptyset$ to $S$ in $\fs (M)$.
\begin{pro}\label{pro:fslat}
If the matroid $M(S)$ has primary flag $T_0\subset\cdots\subset T_k$, 
then the
lattice $\fs (M)$ of free separators of $M$ is equal to the
union of intervals
$\bigcup_{i=1}^k [T_{i-1},T_i]\sub\cf$, where each interval
$[T_{i-1},T_i]\sub\cf$
is a Boolean algebra, given by
$$
[T_{i-1},T_i]\sub\cf\=
\begin{cases}
[T_{i-1},T_i]  & \text{if $T_i$ covers $T_{i-1}$ in $\dlat (M)$,}\\
\{T_{i-1},T_i\} & \text{otherwise,}
\end{cases}
$$
for $1\leq i\leq k$.
\end{pro}
\begin{proof}
By definition, free separators of $M$ are comparable to all
cyclic flats of $M$ and hence comparable to all elements of $\dlat (M)$.
Every free separator is thus contained in one of the intervals
$[T_{i-1},T_i]\sub\cf$, and so $\fs (M) =  
\bigcup_{i=1}^k [T_{i-1},T_i]\sub\cf$.

Suppose that $T_i$ covers $T_{i-1}$ in $\dlat (M)$.  Since 
$T_{i-1}$ and $T_i$ are consecutive pinchpoints of $\dlat (M)$,
and $\dlat (M)$ contains all cyclic flats of $M$, it follows that
any $A\subseteq S$ with $T_{i-1}\subseteq A\subseteq T_i$ is
a free separator.  Hence $[T_{i-1},T_i]\sub\cf =[T_{i-1},T_i]$.

Now suppose that $T_i$ does not cover $T_{i-1}$ in $\dlat (M)$.
Choose some $D\in\dlat (M)$ such that $T_{i-1}\subset D\subset T_i$, and
let $A\in [T_{i-1},T_i]\sub\cf$.  Since $A$ is a free
separator, $A$ must be comparable to $D$.  If $A\subseteq D$,
then the set $\{ E\in\dlat (M)\co A\subseteq E\subset T_i\}$
is nonempty, and thus the intersection $F$ of all elements of this
set is a pinchpoint of $\dlat (M)$ satisfying
$A\subseteq F\subset T_i$.  Since $T_{i-1}$ and $T_i$ are
consecutive pinchpoints of $\dlat (M)$, we therefore
have $A=F=T_{i-1}$.  Similarly, if $D\subseteq A$, it follows
that $A=T_i$.  Hence $[T_{i-1},T_i]\sub\cf =\{T_{i-1},T_i\}$.
\end{proof}

Proposition \ref{pro:fslat} shows, in particular, that $\fs (M)$ is a
sublattice of the Boolean algebra $2^S$, and therefore is a
distributive lattice.  Observe that the first statement of Proposition
\ref{pro:fslat} means that, in addition to being the chain of
pinchpoints in $\dlat (M)$, the primary flag $\calt\sub M$ is also the
chain of all pinchpoints in $\fs (M)$, and the second statement
implies that $\dlat (M)\cap\fs (M)=\calt\sub M$.  If a matroid $M$ has
primary flag $T_0\subset\cdots\subset T_k$, we refer to the minors
$M(T_{i-1},T_i)$ as the {\it primary factors} of $M$, and refer to the
factorization $M=M(T_0,T_1)\frp\cdots\frp M(T_{k-1},T_k)$ as the {\it
  primary factorization} of $M$. 
\begin{thm}\label{thm:setunif}
The sequence of primary factors of a matroid $M$ is the unique sequence 
$M_1,\dots, M_k$ of nonempty matroids such that $M=M_1\frp\cdots\frp M_k$,
each $M_i$ is either irreducible or uniform, and no free product
of consecutive $M_i$'s uniform.
\end{thm}
\begin{proof}
  Suppose that $M(S)$ factors as $M=M_1\frp\cdots\frp M_\ell$.
Let $\cu=\{U_0\subset\cdots\subset U_\ell\}$ be the corresponding
chain in $\fs (M)$, determined by $M_i=M(U_{i-1},U_i)$, 
for $1\leq i\leq\ell$, and let $\calt\sub M=\{T_0\subset\cdots\subset T_k\}$
be the primary flag of $M$. We show that the sequence $M_1,\dots, M_\ell$
has the properties described in the theorem if and only if 
$\cu =\calt\sub M$.

Suppose that $\cu=\calt\sub M$.  By Lemma \ref{lem:fspint}
we have $\fs (M_i)= \fs (M(T_{i-1},T_i))\isom [T_{i-1},T_i]\sub\cf$,
for $1\leq i\leq k$.  If $T_i$ covers $T_{i-1}$ in $\dlat (M)$,
it follows from Proposition \ref{pro:fslat} and Lemma 
\ref{lem:fspuni} that $M_i$ is uniform; and if $T_i$ does not
cover $T_{i-1}$ in $\dlat (M)$, then Proposition \ref{pro:fslat}
and Theorem \ref{thm:irred} imply that $M_i$ is irreducible. 
For $1\leq i\leq k-1$,
we have $M_i\frp M_{i+1}$ = $M(T_{i-1},T_i)\frp M(T_i,T_{i+1})=
M(T_{i-1},T_{i+1})$, and so 
$\fs (M_i\frp M_{i+1})\isom [T_{i-1},T_{i+1}]\sub\cf$,
by Lemma \ref{lem:fspint}.  This interval has a nontrivial 
pinchpoint (namely,
$T_i$), and so is not a Boolean algebra; hence by Lemma \ref{lem:fspuni},
$M_i\frp M_{i+1}$ is not uniform.

For the converse, first note that, since any free separator of $M$
is comparable with all the $T_i$'s, it follows that the union
$\cu\cup\calt\sub M$ is a chain in $\fs (M)$.  Hence if $\calt\not\subseteq
\cu$, we can find $i$ and $j$ such that
$T_j\in [U_{i-1},U_i]\sub\cf$, with $T_j$ not equal to $U_{i-1}$ or
$U_i$.  Then $T_j$ is a nontrivial pinchpoint
of $[U_{i-1},U_i]\sub\cf\isom\fs (M(U_{i-1},U_i))$, and hence
$M_i=M(U_{i-1},U_i)$ is neither uniform nor irreducible.  

Now suppose that $\calt$ is a proper subset of $\cu$.
We can then find some $i$ and $j$
such that $U_j\in [T_{i-1},T_i]\sub\cf$, with $U_j$ not equal to
$T_{i-1}$ or $T_i$.  By Proposition \ref{pro:fslat}, we know that 
$T_i$ covers $T_{i-1}$ in $\dlat (M)$,
from which it follows that $M(T_{i-1},T_i)$ is uniform.  Since
$\calt\subseteq\cu$, we have $T_{i-1}\subseteq U_{j-1}$
and $U_{j+1}\subseteq T_i$; hence the free product 
$M_j\frp M_{j+1}=M(U_{j-1},U_j)\frp 
M(U_j,U_{j+1})=M(U_{j-1},U_{j+1})$ is a minor of $M(T_{i-1},T_i)$
and is thus uniform.
\end{proof}

Theorem \ref{thm:setunif} shows that matroids factor uniquely
as free products of minors that are either irreducible or
``maximally'' uniform. We now wish to consider
factorization of matroids into irreducibles.  
Clearly, given a factorization  
$M(S)=M(U_0,U_1)\frp\cdots\frp M(U_{k-1},U_k)$,
the factors $M(U_{i-1},U_i)$ are all irreducible if and only if
$U_0\subset\cdots\subset U_k$ is a maximal chain in the lattice
of free separators $\fs (M)$.  
If $M(S)=U_{r,n}$ is uniform of rank $r$, then any maximal chain
in $\fs (M)=2^S$, or equivalently, any ordering $s_1,\dots,s_n$ 
of the elements of $S$, gives a factorization 
$$
M\= I(s_1)\frp\cdots\frp I(s_r)\frp Z(s_{r+1})\frp\cdots\frp Z(s_n).
$$
of $M$ into irreducibles (see Example \ref{exa:unif}).  The
factorization of a uniform matroid into irreducibles is thus in general far
from unique.  Up to isomorphism, or course, we do have the unique 
factorization $U_{r,n}=I^r\frp Z^{n-r}$. 
In the next theorem we show that, up to isomorphism,
arbitrary matroids factor uniquely into irreducibles.
\begin{thm}\label{thm:uf}
  If $M\isom M_1\frp\cdots\frp M_k\isom N_1\frp\cdots\frp N_r$,
where all the $M_i$ and $N_j$ are irreducible,
then $k=r$ and $M_i\isom N_j$, for $1\leq i\leq k$. 
\end{thm}
\begin{proof}
Since the sets $T_i$ belonging to the primary flag $\calt\sub M$
of $M$ are all pinchpoints of $\fs (M)$, it follows that any maximal
chain in $\fs (M)$ is a refinement of $\calt\sub M$. Hence
any factorization of $M$ into irreducibles 
can be obtained by starting with the primary factorization
$M\= M(T_0,T_1)\frp\cdots\frp M(T_{\ell-1},T_\ell)$, then factoring each
$M(T_{i-1},T_i)$ into irreducibles.  Since each
$M(T_{i-1},T_i)$ is either irreducible or uniform,
and uniform matroids factor into irreducibles uniquely 
up to isomorphism, it follows that the factorization
of $M$ into irreducibles is unique up to isomorphism.
\end{proof}

The unique factorization theorem \eqref{thm:uf} provides
a quick proof of the following theorem, which was 
the main result in \cite{crsc:fpm}:
\begin{thm}\label{thm:main}
  Suppose that $M(S)\frp N(T)\isom P(U)\frp Q(V)$, where $|S|=|U|$.
Then $M\isom P$ and $N\isom Q$.  
\end{thm}
\begin{proof}
  Since $M\frp N$ and $P\frp Q$ have, up to isomorphism,
the same factorization
into irreducibles, it follows from the fact that $|S|=|U|$ and
$|T|=|V|$, that $M\isom P$ and $N\isom Q$.
\end{proof}

For all $n\geq 0$, denote by $m_n$ and $i_n$, respectively,
the number of isomorphism classes of matroids and irreducible
matroids on $n$ elements, and let $M(t)=\sum_{n\geq 0}m_nt^n$ and
$I(t) = \sum_{n\geq 0}i_n t^n$ be the ordinary generating
functions for these numbers. 
For all $r,k,\geq 0$, denote by 
$m_{r,k}$ and $i_{r,k}$, respectively,
the number of isomorphism classes of matroids and irreducible
matroids having rank $r$ and nullity $k$, and let
$M(x,y) =\sum_{r,k\geq 0}m_{r,k}x^ry^k$ and
$I(x,y) =\sum_{r,k\geq 0}i_{r,k}x^ry^k$.
\begin{cor}\label{cor:gf}
  The generating functions $M(t)$ and $I(t$), and $M(x,y)$
and $I(x,y)$ satisfy
$$
M(t)\=\frac{1}{1- I(t)}\spandsp M(x,y)\=\frac{1}{1- I(x,y)}.
$$
\end{cor}
\begin{proof}
  Unique factorization implies that, for all $n\geq 0$,
$$
m_n\=\sum_{j\geq 0}\!\!\!\!\sum_{\;\;\;\;n_1+\cdots +n_j = n}
\!\!\!i_{n_1}\cdots i_{n_j},
$$
which is the coefficient of $t^n$ in $\sum_{j\geq 0}I(t)^j=
1/(1-I(t))$.  The second equation is proved similarly.
\end{proof}

Using Corollary \ref{cor:gf}, we compute the numbers
$i_n$ and $i_{r,k}$ in terms of the values of 
$m_n$ and $m_{r,k}$, for $n,r+k\leq 8$.  The results are
shown in Tables \ref{table:size} and \ref{table:rn}.

\begin{table}[ht]
\caption{The numbers of nonisomorphic matroids, irreducible matroids, of
size $n$, for $0\leq n\leq 8$:}
\label{table:size}
\begin{tabular}{|c|ccccccccc|}\hline

n & \,0 & 1 & 2 & 3 & 4 & 5 & 6 & 7 & \raisebox{0ex}[2.4ex][1.0ex]{8}  \\
\hline
 matroids & 1 & 2 & 4 & 8 & 17 & 38 & 98 & \!306  & \!\!\!\! 
\raisebox{0ex}[2.35ex][0ex]{1724}   \\
 irreducible matroids & 0 & 2 & 0 & 0 & 1 & 2 & 14 & 66  & 
\raisebox{0ex}[2.15ex][0ex]{891}    \\\hline
 \end{tabular}
\end{table}
\begin{table}[h]
\caption{The numbers of nonisomorphic matroids (left),
irreducible matroids (right), of
rank $r$ and nullity $k$, for $0\leq r+k\leq 8$:}
\label{table:rn}
\begin{tabular}{|c|cccccccclrcccccccc|}\hline
$\!\!\raisebox{-.5ex}{\,r\!}\diagdown\raisebox{.5ex}{\!k}$\!
& 0 & 1 & 2 & 3 & 4 & 5 & \!6 & 7 & \raisebox{0ex}[2.6ex][1.2ex]{8}\,\, 
& \,\,\,\,0 & 1\, & 2 & 3\! & 4 & \!5 & 6\, & 7 & 8 \\
\hline
0 & 1 & 1 & 1 & 1 & 1 & 1 & 1 &  1 & \raisebox{0ex}[2.35ex][0ex]{1}   
& 0 & 1 & 0 & 0 & 0 & 0 & 0 & 0 &  0 \\
1 & 1 & 2 & 3 & 4 & 5 & 6 &  7 &  8 &
& 1 & 0 & 0 & 0 & 0 & 0 & 0  & 0 & \\
2 & 1 & 3 & 7 &13 &23 & 37 & 58 &&&
0 & 0 & 1 & 1 & 3 & 3  & 6  && \\
3 & 1 & 4 &13 &38 &\!\!108 &\!\!325 &&&&
0 & 0 & 1 & 8 & \!30 & \!\!125 &&&  \\
4 & 1 & 5 &23 &\!\!108 &\!\!940 &&&&&
0 & 0 & 3 & \!30 & \!\!629  &&&& \\
5 & 1 & 6 & 37 &\!\!325 &&&&&&
0 & 0 & 3  & \!\!125  &&&&& \\
6 & 1 &  7 & 58 &&&&&&&
0 & 0  & 6  &&&&&& \\
7 & 1 &  8 &&&&&&&&
0  & 0  &&&&&&& \\
8 & 1 &&&&&&&&&
0  &&&&&&&& \\
\hline
\end{tabular}
\vspace{2ex}
\end{table}

The two matroids of size one, namely, the point \pt\ and loop \lp,
are irreducible, and no matroid of size two or three
is irreducible.  The unique irreducible matroid on four elements is
the pair of double points
$\un 12\oplus\un 12$.  The two irreducible matroids on five elements
are $\un 13\oplus\un 12$ and its dual $\un 23\oplus\un 12$.
On six elements, the irreducibles of rank two are
$\un 14\oplus\un 12$, $\un 13\oplus\un 13$ and the truncation
$\trun{}{(\un 12\oplus\un 12\oplus\un 12)}$, which consists 
of three collinear
double points.  The duals of these matroids, 
$\un 34\oplus\un 12$, $\un 23\oplus\un 23$ and $\lift{}{(\un 12\oplus\un 12
\oplus\un 12)}$, are the six-element irreducibles of rank four.
Finally, on six elements in rank three, the irreducible matroids consist
of $\un 24\oplus\un 12$, $\un 12\oplus\un 12\oplus\un 12$, $\un 13\oplus
\un 23$, and $\un 23'\oplus \un 12$, where $\un 23'$ is the three-point
line $\un 23$, with one point doubled, together with the four matroids
shown below:

\vspace{.9em}
\mbox
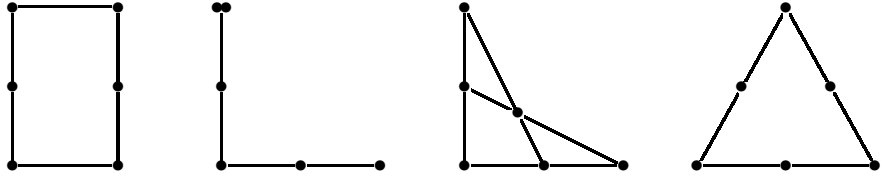
{
\begin{diagram}[PostScript=Rokicki,abut,height=1.0em,width=1.0em,tight,thick]
\bullet&\rLine&&&\bullet\\ \\
\dLine&&&&\dLine \\
\bullet&&&&\bullet\\ \\
\dLine&&&&\dLine \\
\bullet&\rLine&&&\bullet
\end{diagram}}
\mbox{
\begin{diagram}[PostScript=Rokicki,abut,height=1.0em,width=1.0em,tight,thick]
&&&\mbox{$\bullet\!\bullet$}\\ \\
&&&\dLine \\
&&&\bullet\\ \\
&&&\dLine \\
&&&\bullet&\rLine&&\bullet&\rLine&&\bullet
\end{diagram}}\;
\mbox{
\begin{diagram}[PostScript=Rokicki,abut,height=1.0em,width=1em,tight,thick]
&&\bullet\\ 
&&&\rdLine(2,4)\\ 
&&\dLine \\
&&\bullet\\
&&&\rdLine(2,1)&\bullet \\ 
&&\dLine&&&\rdLine (4,2)\rdLine(1,2)\\
&&\bullet&\rLine&&\bullet&\rLine&&\bullet
\end{diagram}}\;\,
\mbox{
\begin{diagram}[PostScript=Rokicki,abut,height=1.5em,width=1.68em,tight,thick]
&&&\bullet\\ 
&&\ldLine(1,2)&&\rdLine(1,2)\\ 
&&\bullet &&\bullet\\ 
&\ldLine(1,2)&&&&\rdLine(1,2)\\ 
&\bullet&\rLine&\bullet&\rLine&\bullet
\end{diagram}}

\vspace{2ex}
\noindent
Since the dual of an irreducible matroid is irreducible,
the set of rank-three irreducible matroids on six elements must be closed
under duality; in fact, each matroid in this set is self-dual.

\section{The minor coalgebra}\label{sec:algebra}

In this section, and the next,
we work over some commutative ring \K\ with unit.  All
modules, algebras and coalgebras are over \K, all maps between such
objects are assumed to be \K-linear, and all tensor products are taken
over \K.  Given a family of matroids \cm, we denote by \frmod\cm\ the
free \K-module having as basis all isomorphism classes of matroids
belonging to \cm.  In what follows, we shall not distinguish
notationally between a matroid $M$ and its isomorphism class, or
between a family of matroids \cm\ and the set of isomorphism classes
of matroids belonging to \cm; it should always be clear from the
context which is meant.

If \cm\ is a minor-closed family, then the {\it minor
coalgebra} (\cite{sc:iha}, \cite{crsc:fsa}) of \cm\ is
the free module \frmod\cm, equipped with {\it restriction-contraction
coproduct\/} $\delta$ determined by
$$
\delta (M)\= \sum\sub{A\subseteq S}M|A\o M/A,
$$
and counit determined by $\epsilon (M) = \delta\sub{\emptyset, S}$,
for all $M=M(S)$ in $\cm$. 
If \cm\ is also closed under formation of
direct sums, then \frmod\cm\ is a Hopf algebra, with product determined
on the basis \cm\ by direct sum.  
For any minor-closed family \cm, the coalgebra \frmod\cm\ is bigraded,
with homogeneous component $\frmod\cm_{r,k}$ spanned by all
isomorphism classes of matroids in \cm\ having rank $r$ and nullity
$k$.  When \cm\ is also closed under direct sum, this is a Hopf algebra  
bigrading.

For all matroids $N_1$, $N_2$ and $M=M(S)$, the {\it section coefficient} 
\sect M{N_1}{N_2}\ is
the number of subsets $A$ of $S$ such that $M|A\isom N_1$ and
$M/A\isom N_2$; hence if \cm\ is a minor-closed family, the 
restriction-contraction coproduct satisfies
\begin{equation}
  \label{eq:seccop}
\delta (M)\=\sum\sub{N_1,N_2}\sect M{N_1}{N_2} N_1\o N_2,  
\end{equation}
 for all $M\in\cm$,
where the sum is taken over all (isomorphism classes of) matroids
 $N_1$ and $N_2$.  More generally, for matroids $N_1,\dots, N_k$ and 
$M=M(S)$, the {\it multisection coefficient} \msect M{N_1}{N_k} is
the number of sequences $(S_0,\dots,S_k)$ such that
 $\emptyset =S_0\subseteq\cdots
\subseteq S_k=S$ and the minor $M(S_{i-1},S_i)$ is isomorphic to 
$N_i$, for $1\leq i\leq k$.
The iterated coproduct 
$\delta^{k-1}\co\frmod\cm\rta\frmod\cm\o\cdots\o\frmod\cm$
is thus determined by
$$
\delta^{k-1} (M)\=\sum\sub{N_1,\dots,N_k}\msect M{N_1}{N_k}
N_1\o\cdots\o N_k,
$$
for all $M\in\cm$.

For any family of matroids \cm, we define a pairing $\pr\cdot\cdot\co
\frmod\cm\times\frmod\cm\rta\K$ by setting $\pr MN=\delta\sub{M,N}$, 
for all $M,N\in\cm$, and
thus identify the graded dual module $\frmod\cm\dual$ with
the free module \frmod\cm.  In the case that \cm\ is
minor-closed, we refer to the (graded) dual algebra 
$\frmod\cm\dual$ as the {\it minor algebra} of \cm; the product
in the minor algebra is thus determined by
$$
M \cdot N\=\sum\sub{L\in\cm}\sect LMN \,L,
$$
for all $M,N\in\cm$.

We partially order the set of all isomorphism classes of matroids
by setting $M\geq N$ if and only if there exists a bijective
weak map from $M$ to $N$.  The following result provides us
with critical necessary conditions for a matroid to 
appear in a given product of matroids 
in $\frmod\cm\dual$. 
\begin{pro}\label{pro:initial}
For all matroids $L$,$M$ and $N$,
$$
\sect LMN\neq 0 \quad\Longrightarrow\quad 
M\oplus N\,\leq\, L\,\leq\, M\frp N.
$$
\end{pro}
\begin{proof}
   Suppose that $M=M(S)$ and $N=N(T)$. Given a matroid $L$ such that
$\sect LMN\neq 0$ we may assume that $L=L(S+T)$, where $L|S=M$ and $L/S=N$.  
The semimodularity of $\rho\sub L$ implies that 
$\rkl\as +\rkl{S\cup A}\leq\rkl S +\rkl A$, for all $A\subseteq S+T$, 
and so $\rho\sub{M\oplus N}(A) =\rkm\as +\rkn\at 
=\rkl\as +\rkl{S\cup A}-\rkl S
\leq\rkl A$, and hence the identity on $S+T$ is a weak map
 $L\rta M\oplus N$.   On the other hand, according to
Proposition \ref{pro:univ}, the identity on $S+T$ is a weak map
$M\frp N\rta L$; hence $M\oplus N\leq L\leq M\frp N$.
\end{proof}

\noindent
Similarly, using Proposition \ref{pro:guniv} instead of
Proposition \ref{pro:univ}, we obtain 
\begin{equation}
  \label{eq:initial}
  \msect L{M_1}{M_k}\neq 0 \quad\Longrightarrow\quad
M_1\oplus\cdots\oplus M_k \,\leq\, L\,\leq\, M_1\frp\cdots\frp M_k,
\end{equation}
for all $L$ and $M_1,\dots, M_k\in\cm$.

The following example shows that the converse of Proposition
\ref{pro:initial} does not hold.
\begin{exa}
  Suppose $L$ is the rank $4$ matroid on the set
$U=\{a,b,c,d,e,f,g\}$ pictured below.
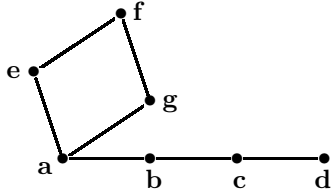
\newcommand{\exlab}[1]
{\hspace{2.ex}\text{\raisebox{-.8ex}{$\bf #1$}}}
  \begin{diagram}[PostScript=Rokicki,abut,height=1.1em,width=1.1em,tight,thick,balance]
&&&&\bullet&{\!\!\!\!\!\bf f}\\
&&&\ruLine(3,2)&&\rdLine(1,3)\\
{\bf e}\!\!\!\!&\bullet\\
&&\rdLine(1,3)&&&\bullet&\!\!\!\!{\bf g}\\
&&&&\ruLine(3,2)\\
&\exlab{a}&\bullet&\rLine&&\bullet&\rLine&&\bullet&\rLine&&\bullet&
{\hspace{-20.2ex}\text{\raisebox{-2.ex}
{$\bf b$\hspace{6.2ex}$\bf c$\hspace{6.ex}$\bf d$}}}
\end{diagram}

\vspace{1.ex}
\noindent
If $M$ is a three point line on the set $\{a,b,c\}$, and
$N$ is a four point line on $\{d,e,f,g\}$, then the free
product $M\frp N$ consists of a three point line on $\{a,b,c\}$,
together with points  $d$, $e$, $f$, $g$ in general position
in $3$-space, and the identity map on $U$ is thus 
a weak map $M\frp N\rta L$.  Now if $M'$ a three point line 
on $\{e,f,g\}$ and $N'$ is a four point line on $\{a,b,c,d\}$, 
then the identity on $U$
is a weak map $L\rta M'\oplus N'$.  Since $M\isom M'$ and $N\isom N'$,
we thus have $M\oplus N\leq L\leq M\frp N$.  But $L$ has no
three point line as a restriction with a four point line as
complementary contraction, and so $\sect LMN=0$.
\end{exa}

If a family \cm\ is closed under formation of free products
then  \frmod\cm, with product determined 
by the free product on the basis \cm, is an associative algebra. 
We denote \frmod\cm, equipped with this algebra structure, by
$\frmod\cm\sub\frp$.  
\begin{pro}\label{pro:free}
  If \cm\ is a free product-closed family of matroids, then
the algebra $\frmod\cm\sub\frp$ is free, generated by the set of
irreducible matroids belonging to \cm.
\end{pro}
\begin{proof}
Because the set \cm\ is a basis
for $\frmod\cm\sub\frp$, the result follows directly from 
unique factorization, Theorem \ref{thm:uf}.
\end{proof}

For all matroids $M$ and $N$, we denote by $c(N,M)$ the section
coefficient \msect N{M_1}{M_k}, where $M_1,\dots, M_k$ is the 
sequence of irreducible factors
of $M$. 
\begin{thm}\label{thm:free}
 Suppose that \cm\ is a family of matroids that is closed under
formation of minors and free products. If \K\ is a field of characteristic 
zero, then the minor algebra $\frmod\cm\dual$ is free, generated by the 
set of irreducible matroids belonging to \cm.
\end{thm}
\begin{proof}
  For each matroid $M$ belonging to \cm, let $P\sub M$ denote
the product $M_1\cdots M_k$ in $\frmod\cm\dual$,
where $M_1,\dots, M_k$ is the sequence of irreducible factors
of $M$.  We can write 
$$
P\sub M\=\sum\sub{N} c(N,M)\, N,
$$
where, by \eqref{eq:initial}, the sum is taken over all $N\in\cm$ 
such that $N\leq M$ in the weak order.  Since $c(M,M)\neq 0$, 
for all matroids $M$, and \K\ is a field of characteristic zero,
it thus follows from the fact that \cm\ is a basis for
$\frmod\cm\dual$ that $\{P\sub M\,\co M\in\cm\}$ 
is also a basis for $\frmod\cm\dual$.  The map $\frmod\cm\sub\frp\rta
\frmod\cm\dual$ determined by $M\mapsto P\sub M$, which is clearly an algebra
homomorphism, is thus bijective and hence an algebra isomorphism.
Since $P\sub M=M$, whenever $M\in\cm$ is irreducible, the result follows from
Proposition \ref{pro:free}.
\end{proof}
\begin{exa}
  The family \cm\ of all matroids is minor-closed and
closed under free product.  Hence $\frmod\cm\dual$ is the free
algebra generated by the set of all (isomorphism classes of) irreducible 
matroids.  
\end{exa}
\begin{exa}
  The family \cf\ of freedom matroids (see Example \ref{exa:freedom}) is
minor-closed and closed under free product.  Since all
freedom matroids can be expressed as free products of
points and  loops, it follows that $\frmod\cf\dual$ is
the free algebra generated by \pt\ and \lp.
\end{exa}
\begin{exa}
For any field $F$, the class $\cm\sub F$ of all $F$-representable matroids
is minor-closed.  It follows from Proposition \ref{pro:repprod}
that if $F$ is infinite then $\cm\sub F$ is also closed under formation
of free products.
\end{exa}
\begin{exa}
It follows from Proposition \ref{pro:trans} that the family
\calt\ of all transversal matroids is closed under
formation of free products.  However, since contractions of transversal
matroids are not in general transversal, \calt\ is not minor-closed.
\end{exa}
\begin{pro}
  If a family \cm\ of matroids is minor-closed and closed
under formation of free products, then \cm\ is also closed
under the lift and truncation operations.
\end{pro}
\begin{proof}
  Suppose that \cm\ is minor-closed and closed under formation
of free products.  If \cm\ is the class of all free matroids
or the class of all zero matroids, or consists only of the empty
matroid, then \cm\ is closed under lift and truncation.  If
\cm\ is none of the above, then it must contain the matroids
\pt\ and \lp.  By Proposition \ref{pro:delcon},
we have
$$
\lift{}M\= (\pt\frp M(S))|S \sands \trun{}N\= (M\frp Z(a))/a,
$$
for any matroid $M=M(S)$.  Hence, if $M$ belongs to \cm\, then
so do $\lift{}M$ and $\trun{}M$.
\end{proof}

Suppose that \cm\ and \K\ satisfy the hypotheses of Theorem \ref
{thm:free}, and that \cm\ is partially ordered by the weak order.
The fact that $c(M,N)\neq 0$ implies $M\leq N$, for all $M,N\in\cm$, 
means that we may regard $c$ as an element of the incidence
algebra $I(\cm)$ of the poset \cm.  Since $c(M,M)$ is invertible in \K,
for all $M$, it follows that $c$ is invertible in $I(\cm)$, the
inverse given recursively by $c^{-1}(M,M)=
c(M,M)^{-1}$, for $M\in\cm$, and 
$$
 c^{-1}(M,N)\= -\,c(N,N)^{-1}\!\!\!\sum\sub{M\leq P<N}c^{-1}(M,P)\,c(P,N),
$$
for all $M<N$ in \cm.
The inverse of the change of basis map
$M\mapsto P\sub M$ is thus given by
$$
M \=\sum\sub{N\leq M} c^{-1}(N,M)P\sub N,
$$
for all $M\in\cm$.
Let $\{Q\sub M\co M\in\cm\}$ be the basis of \frmod\cm\ determined
by $\langle Q\sub M, P\sub N\rangle = \delta\sub{M,N}$, for all
$M,N\in\cm$.  Observe that $Q\sub M$ is satisfies
\begin{equation}\label{eq:prim}
Q\sub M\=\sum\sub{N\geq M}c^{-1}(M,N)N,
\end{equation}
for all $M\in\cm$.  Before stating the next theorem, which
is dual to Theorem \ref{thm:free}, we note that, for any
minor-closed family \cm,  the minor
coalgebra \frmod\cm\ is connected, with the empty matroid as
unique group-like element.  In particular, it follows that the notion
of primitive element of \frmod\cm\ is unambiguous.
\begin{thm}
  Suppose that \cm\ is a family of matroids that is closed
under formation of minors and free products.  If \K\ is a field
of characteristic zero, then the minor coalgebra \frmod\cm\ is cofree.  
The set $\{Q\sub M\co\text{$M\in\cm$ is irreducible}\}$ is
a basis for the subspace of primitive elements of \frmod\cm.
\end{thm}
\begin{proof}
The fact that \frmod\cm\ is cofree is equivalent to the fact
that \frmod\cm\dual\ is free, which was shown in Theorem \ref{thm:free}.
Let $\varphi\co\frmod\cm\sub\frp\rta\frmod\cm\dual$ be the algebra
isomorphism used in the proof of Theorem \ref{thm:free}, given 
by $M\mapsto P\sub M$, for all $M\in\cm$.  The
transpose $\varphi\dual\co\frmod\cm\rta\frmod\cm\sub\frp\dual$
is thus a coalgebra isomorphism. For all $M,N\in\cm$, we have
$\langle \varphi\dual (Q\sub M), N\rangle = 
\langle Q\sub M, \varphi (N)\rangle =
\langle Q\sub M, P\sub N\rangle = \delta\sub{M,N}$, and hence
$\varphi\dual (Q\sub M) = M$.  Since the
set of all irreducible $M\in\cm$ is a basis for the subspace of
primitive elements of $\frmod\cm\sub\frp\dual$, it follows
that $\{Q\sub M\co\text{$M\in\cm$ is irreducible}\}$ is a basis
for the subspace of primitive elements of $\frmod\cm$.
\end{proof}
\begin{exa}
  Suppose that \cm\ is closed under formation of minors and
free products, and that \cm\ contains the irreducible matroid
$D=U_{1,2}\oplus U_{1,2}$, consisting of two double points.
Since \cm\ is minor-closed, it contains the (irreducible)
single-element matroids $I$ and $Z$.  Since \cm\ is also closed
under free product, it follows from Table \ref{table:size} and
unique factorization that \cm\ contains all matroids of size
less than or equal to four (all such matroids, except $D$, being
free products of $I$'s and $Z$'s).  

   It is clear from Equation \ref{eq:prim} that the primitive
elements $Q\sub I$ and $Q\sub Z$ in \frmod\cm\ are equal
to $I$ and $Z$, respectively.  In order compute $Q\sub D$, we first
observe that $\{ N\co\text{$N > D$ in \cm}\}$ consists of
the two matroids $U_{2,4}=I\frp I\frp Z\frp Z$ and $P=I\frp Z\frp I\frp Z$.
Since $P$ is a three point line, with one point doubled, we have
$D\leq P\leq U_{2,4}$.  The multisection coefficients $c(M,N)$,
for all $M,N\geq D$, are given by the matrix
$$
\bordermatrix{
&{\raisebox{0ex}{$\scriptstyle D$}}&
{\raisebox{0ex}{$\scriptstyle P$}}&
{\raisebox{0ex}{$\scriptstyle U_{2,4}$}}\cr
{\scriptstyle D} &  1& 8&16\cr
{\scriptstyle P} &  0 &4 &20\cr
{\scriptstyle U_{2,4}} &  0 &0 &24\cr
}
$$\\

\noindent
and the numbers $c^{-1}(M,N)$, for $M,N\geq D$, are thus
given by the inverse matrix

$$
\frac1{24}\left(\begin{array}{crr}
24 & -48 & 24\\
0 & 6 & -5\\
0&0&1
\end{array}
\right).
$$\\
\noindent
Hence $Q\sub D= D-2P+U_{2,4}$.
\end{exa}

\section{A new twist}\label{sec:twist}

If a family of matroids \cm\ is both minor and free product-closed,
then the module \frmod\cm\ has both the structure of a
(free) associative algebra, under free product, and a coassociative
coalgebra, with restriction-contraction coproduct.  Moreover,
according to Theorem \ref{thm:free}, when the ring of scalars 
is a field of characteristic zero, these structures
are dual to one another.  In this section we show that
free product and restriction-contraction coproduct are
compatible in the sense that \frmod\cm\ is a Hopf
algebra in an appropriate braided monoidal category.

By a {\it matroid module}, we shall mean
a free module \frmod\cm, where \cm\ is a family of matroids that
is closed under formation of lifts and truncations.
Given matroid modules $V=\frmod\cm$ and $W=\frmod\cn$,
we define the
{\it twist map\/} $\tau=\tau\sub{V,W}\co V\o W\rta W\o V$ by
\begin{equation}
  \label{eq:twist}
  \tau (M\o N)\= \lift{\rho(M)}N\o\trun{\nu(N)}M,
\end{equation}
for all $M\in\cm$ and $N\in\cn$.  If the families \cm\ and \cn\ are
also closed under formation of free products, we use the twist map to
extend the definition of the free product to a binary operation on
$V\o W$:
\begin{equation}
\label{eq:frp}
(M\o N)\frp (P\o Q)\= (M\frp\lift{\rho(N)}P)\o (\trun{\nu(P)}N\frp Q),
\end{equation}
for all $M,P\in\cm$ and $N,Q\in\cn$.
\begin{pro}\label{pro:associative}
  For all families \cm\ and \cn, closed under free product,
lift and truncation, the product $\frp$ given by Equation
\ref{eq:frp} is an associative operation on $\frmod\cm\o\frmod\cn$.
\end{pro}
\begin{proof}
Suppose that $M_i\in\cm$ and $N_i\in\cn$, and let $\nu_i=\nu(M_i)$ and
$\rho_i=\rho(N_i)$, for $1\leq i\leq 3$.  Then 
\begin{align*}
[(M_1\o N_1)\frp (M_2\o N_2)] &\frp (M_3\o N_3)\\
&\=[(M_1\frp\lift{\rho_1}M_2)\o(\trun{\nu_2}N_1\frp N_2)] 
\frp (M_3\o N_3)\\ 
&\= (M_1\frp \lift{\rho_1}M_2\frp\lift iM_3)\o
   (\trun{\nu_3}(\trun{\nu_2}N_1\frp N_2)\frp N_3)\\ 
&\= (M_1\frp\lift{\rho_1}M_2\frp\lift iM_3)\o
   (\trun kN_1\frp\trun{\nu_3}N_2\frp N_3),
\end{align*}
Where $i =\rho (\trun{\nu_2}N_1\frp N_2) =
\rho_2+\max\{\rho_1-\nu_2,\,0\}$ and, by Equation \ref{eq:ltfrp},
we have $k=\nu_2+\max\{\nu_3-\rho_2,\, 0\}$.  
On the other hand,
\begin{align*}
(M_1\o N_1)\frp [(M_2\o N_2)&\frp (M_3\o N_3)]\\
&\= (M_1\o N_1)\frp [(M_2\frp\lift{\rho_2}M_3)\o 
(\trun{\nu_3}N_2\frp N_3)]\\ 
&\= (M_1\frp \lift{\rho_1}(M_2\frp\lift{\rho_2}M_3))\o
(\trun jN_1\frp\trun{\nu_3}N_2\frp N_3)  \\ 
&\= (M_1\frp\lift{\rho_1}M_2\frp\lift sM_3 )\o
(\trun jN_1\frp\trun{\nu_3}N_2)\frp N_3),  
\end{align*}
  where $j=\nu(M_2\frp\lift{\rho_2}M_3) = \nu_2 +\max\{\nu_3-\rho_2,\, 0\}$
and, by Equation \ref{eq:ltfrp}, we have $s=\rho_2+\max\{\rho_1-\nu_2,\,0\}$.
Since $s=i$ and $j=k$, the two parenthesizations of $(M_1\o N_1)
\frp (M_2\o N_2)\frp (M_3\o N_3)$
are thus equal.
\end{proof}

\begin{pro}\label{pro:bialgebra}
 If the family \cm\ is minor and free product-closed (and thus also 
closed under lift and truncation), then the restriction-contraction
coproduct $\delta$ is compatible with the free product on \frmod\cm, in
the sense that $\delta\co\frmod\cm\rta\frmod\cm\o\frmod\cm$ is
an algebra map. 
\end{pro}
\begin{proof}
 Suppose that $M(S)$ and $N(T)$ belong to \cm.  Using Proposition
\ref{pro:delcon}, we compute the coproduct of $M\frp N$:
 \begin{align*}
   \delta (M\frp N) 
&\=\sum\sub{A\subseteq S+T}(M\frp N)|A \o (M\frp N)/A\\
&\=\sum\sub{A\subseteq S+T}(M|\as \frp \lift{\rlm\as}N|\at)
\o (\trun{\nln\at}M/\as\frp N/\at)\\
&\=\sum\sub{A\subseteq S+T}(M|\as \frp \lift{\rho(M/\as)}N|\at)
\o (\trun{\nu(N|\at)}M/\as\frp N/\at)\\
 &\=\sum\sub{A\subseteq S+T}(M|\as \o M/\as)\frp (N|\at\o N/\at),
\end{align*}
which is equal to $\delta (M)\frp\delta (N)$.
\end{proof}

We conclude by outlining a categorical framework for these
results.  Let \mcat\ be the category whose objects are
bigraded \K-modules $V=\bigoplus_{r,k\geq 0}V_{r,k}$, equipped
with linear operators $\li=\li\sub V$ and $\tr=\tr\sub V$ satisfying

\begin{enumerate}[(i)]
\item $\li\co V_{r,k}\rta V_{r+1,k-1}$, if $k>0$,\sands
$\li|V_{r,0}=\idt{V_{r,0}}$,\\
\item $\tr\co V_{r,k}\rta V_{r-1,k+1}$, if $r>0$,\sands
 $\tr|V_{0,k}=\idt{V_{0,k}}$,\\
\item $\tr\li=\li\tr$, when restricted to $\bigoplus_{r,k\geq 1}V_{r,k}$.\\
\end{enumerate}
We assume that each homogenous component $V_{r,k}$ is a free
\K-module of finite rank and that $V_{r,0}$ and $V_{0,k}$
have rank one, for all $r,k\geq 0$.  For homogeneous
$x\in V$, we write $\rho (x)=r$ and $\nu (x)=k$ to indicate
that $x$ belongs to $V_{r,k}$.
The morphisms of \mcat\
are the \K-linear maps which commute with $\li$ and $\tr$.  
For all objects $V$ and $W$ in \mcat, we suppose that
the tensor product $V\o W$ is bigraded in the usual manner, with
$$
(V\o W)_{r,k}\= 
\bigoplus_{\stackrel{r_1+r_2=r}{\scriptscriptstyle k_1+k_2=k}}
(V_{r_1,k_1}\o W_{r_2,k_2}),
$$
for all $r,k\geq 0$, and the operators $\li=\li\sub{V\o W}$
and $\tr=\tr\sub{V\o W}$ satisfy
$$
\lift{}{(x\o y)}\=
\begin{cases}
  (\lift{}x)\o y & \text{if $\nu (x)>0$,}\\
 x\o\lift{}y & \text{if $\nu (x)=0$,}
\end{cases}
$$
and
$$
\trun{}{(x\o y)}\=
\begin{cases}
  x\o \trun{}y & \text{if $\rho (y)>0$,}\\
 (\trun{}x)\o y & \text{if $\rho (y)=0$,}
\end{cases}
$$
for all homogeneous $x\in V$ and $y\in W$; hence \mcat\ is
a monoidal category.  For all objects $V$ and $W$ in \mcat\
we define the  twist map $\tau=\tau\sub{V,W}\co V\o W\rta W\o V$ as
in Equation \ref{eq:twist}, that is, by
$\tau (x\o y)\= \lift{\rho(x)}y\o\trun{\nu(y)}x$,
for homogeneous $x\in V$ and $y\in W$.  
It is readily verified that
the twist maps $\tau\sub{V,W}$ commute with
the operators \lift{}\ and \trun{}, and so are morphisms
in \mcat; furthermore, the maps
$\tau\sub{V,W}$ are the components of
a natural transformation $\tau\co \o\Rightarrow\o^{\rm op}$, that
is,  $(g\o f)\circ\tau\sub{V,W}=\tau\sub{V',W'}\circ (f\o g)$,
for all morphisms $f\co V\rta V'$ and $g\co W\rta W'$ in \mcat.
It is then a simple matter to verify that 
the natural transformation $\tau$ satisfies the braid relations:
$$
\tau\sub{U\o V,W}\= (\tau\sub{U,W}\o 1\sub V)\circ
(1\sub U\o\tau\sub{V,W})\sands
\tau\sub{U, V\o W}\= (1\sub V\o\tau\sub{U,W})\circ
(\tau\sub{U,V}\o 1\sub W),
$$
for all objects $U,V,W$. 
Note that the maps $\tau\sub{V,W}$ are not necessarily isomorphisms
in \mcat\ (because different matroids may have the same lifts
or truncations).   Hence, as long as we allow a notion of
braiding that does not require the component morphisms to be isomorphisms, 
it follows that \mcat\ is a braided monoidal category.

We regard each matroid module \frmod\cm\ as an object of \mcat,
bigraded by rank and nullity, with operators \li\ and \tr\ 
determined by lift and truncation on the basis \cm.  If $V=\frmod\cm$,
and the family of matroids \cm\ is closed under free product, as well
as lift and truncation, then it follows immediately from 
Proposition \ref{pro:ltfrp} and the definition of \li\ and \tr\
on $V\o V$ that the map $\mu\sub V\co V\o V\rta V$ given by
$M\o N\mapsto M\frp N$, for all $M,N\in\cm$, is a morphism
in \mcat, and hence $V$ is a monoid object in \mcat.

Suppose that $V=\frmod\cm$ and $W=\frmod\cn$ are matroid
modules with $\cm$ and $\cn$ free product-closed. 
The operation $\frp$ on $V\o W$ defined by Equation
\ref{eq:frp} is the composition 
$\mu\sub{V\o W}
=(\mu\sub V\o\mu\sub W)\circ(1\sub V\o\tau\sub{V,W}\o 1\sub W)$,
 which is the standard monoid structure on
the product of monoid objects in a braided monoidal
category.  Associativity of $\mu\sub{V\o W}$ (our Proposition
\ref{pro:associative}) follows immediately from the braid
relations and the associativity of $\mu\sub V$ and $\mu\sub W$.

Finally, we note that if $V=\frmod\cm$ is a matroid module, where
\cm\ is minor-closed, then the restriction-contraction coproduct
$\delta\co V\rta V\o V$ commutes with \li\ and \tr, and so $V$ is
a comonoid object in \mcat.  If \cm\ is also closed under free
product, then Proposition \ref{pro:bialgebra} says that
$V$ is a bialgebra in the braided monoidal category \mcat.
Since $V$ is a connected bialgebra, it is in fact
a Hopf algebra, with antipode given by the usual formula.
Furthermore, it follows from the proof of Theorem \ref{thm:free}
that this Hopf algebra is self-dual.

\providecommand{\bysame}{\leavevmode\hbox to3em{\hrulefill}\thinspace}
\providecommand{\MR}{\relax\ifhmode\unskip\space\fi MR }
\providecommand{\MRhref}[2]{%
  \href{http://www.ams.org/mathscinet-getitem?mr=#1}{#2}
}
\providecommand{\href}[2]{#2}


\begin{thebibliography}{10}

\bibitem{ar:tcm}
Frederico Ardila, \emph{The {C}atalan matroid}, Journal of Combinatorial 
Theory A \textbf{104} (2003), 49--62.

\bibitem{bodeno:lpm}
Joseph Bonin, Anna de~Mier, and Marc Noy, \emph{Lattice path matroids:
enumerative aspects and {T}utte polynomials}, Journal of Combinatorial 
Theory A \textbf{104} (2003), 63--94.

\bibitem{cr:see}
Henry Crapo, \emph{Single-element extensions of matroids}, Journal of 
research, National Bureau of Standards \textbf{69B} (1965), 55--66.

\bibitem{crsc:fpm}
Henry Crapo and William Schmitt, \emph{The free product of matroids}, 
accepted for publication in the European Journal of Combinatorics (2004), 
arXiv:math.CO/0409080.

\bibitem{crsc:fsa}
\bysame, \emph{A free subalgebra of the algebra of matroids}, accepted for
publication in the European Journal of Combinatorics (2004), arXiv:math.CO/0409028.

\bibitem{hi:smg}
Denis Higgs, \emph{Strong maps of geometries}, Journal of Combinatrial Theory
\textbf{5} (1968), 185--191.

\bibitem{ox:mt}
James Oxley, \emph{Matroid {T}heory}, Oxford University Press, Oxford, 1992.

\bibitem{sc:iha}
William Schmitt, \emph{Incidence {H}opf algebras}, Journal of Pure and 
Applied Algebra \textbf{96} (1994), 299--330.

\bibitem{we:bnm}
Dominic J.~A. Welsh, \emph{A bound for the number of matroids}, Journal of
Combinatorial Theory \textbf{6} (1969), 313--316.

\bibitem{we:mt}
\bysame, \emph{Matroid {T}heory}, Academic Press, London, 1976.

\end{thebibliography}
\end{document}